\documentstyle[12pt,epsf,amssymb]{article}
\newtheorem{proposition}{Proposition}
\title{Poisson structure on moduli of flat connections \\
on Riemann surfaces and $r$-matrix}
\author{V.V.Fock and A.A.Rosly \\
 {\em Institute of Theoretical and Experimental Physics,}
\\ {\em B.Cheremushkinskaya 25,} \\ {\em 117259,  Moscow, Russia}
\date{}}
\begin{document}
\newcommand{\noeq}[1]{\refstepcounter{equation}\eqno(#1\arabic{equation})}
\maketitle
\newcommand{\bird}[1]{#1^\vee}
\begin{abstract}
We consider the space of graph connections (lattice gauge fields) which
can be endowed with
a Poisson structure in terms of a {\em ciliated fat graph}. (A ciliated
fat graph is a graph with a fixed linear order of ends of  edges at  each
vertex.) Our aim is however to study the Poisson structure on the
moduli space of locally flat vector bundles on a Riemann surface with
holes (i.e. with boundary).
It is shown that this moduli space can be obtained as a
quotient of the space of graph connections
by the Poisson action of a lattice gauge group endowed with a
Poisson-Lie structure.
The present paper contains
as a part an updated version of a 1992 preprint
\cite{FR} which we decided still deserves publishing. We have removed some
obsolete inessential remarks and added some newer ones.  \end{abstract}

\newpage
\section{Introduction}

The  moduli space of flat $G$-bundles on a Riemann surface is the classical
phase space for the Chern-Simons gauge theory and, thus, it  is  in  a  sense
the classical limit of the WZNW conformal field theory. This means that
quantizing it one can get a space of quantum states which turns out to be
isomorphic  to the space of conformal blocks of the corresponding WZNW theory.
This statement has been checked by several authors (cf. \cite{ADPW},\cite{Wit})
with help of different quantization methods. On the other hand, the moduli
spaces of flat  bundles as well as closely related to them moduli spaces of
holomorphic bundles (cf. \cite{NR}) attracted much attention from a purely
mathematical point of view (cf. \cite{AB,H}).

   In section  \ref{ps} we discuss in
detail the canonical Poisson structure on the moduli space of flat bundles
on Riemann surfaces with holes (i.e. with boundary). In section \ref{fdm}
we construct a Poisson structure on the space of  graph connections in such a
way that the action of the graph gauge group is Poissonian  with respect to
an appropriate nontrivial Poisson-Lie structure. The considerations of this
section are inspired mainly by constructions of  refs. \cite{STS,AFSV,AFS}
where a discrete analogue of current algebra was suggested  and investigated.
Then we prove that the quotient of the space of graph connections  by the
gauge group coincides with the moduli space of flat connections  on  a
Riemann surface determined by the graph.

One of the main aims of our preprint \cite{FR} (see also ref. \cite{FR2})
was to give a description of
the moduli space of flat connections in the form ready for quantisation.
We are not going to discuss this problem here for since that time
some progress towards this direction was made,
see refs. \cite{AGSch,BuffRo,Frol}. The interested reader
can find details and discussions in those papers.  Note only that the
result of quantisation is a noncommutative algebra having  WZNW  model
conformal blocks as its representation space and being functorial with respect
to the imbeddings of surfaces.

For the readers interested in a more understandable and detailed presentation
we would like to recommend a very good review by M.Audin
\cite{Audin}\footnote{It worths to look through ref. \cite{Audin} not only
because of a very transparent presentation, but also because of very nice
pictures there.}.

\section{Ciliated fat graphs and Poisson manifolds} \label{cfg}

The moduli space of flat connections on a compact Riemann surface
is by definition a subquotient of a
topologically trivial space of all connections. This description is
useful also
since a nontrivial Poisson manifold  (which  is  the  moduli
space, or an orbifold, to be more precise) is
represented as a result of a reduction of a  trivial
symplectic manifold (see sect. \ref{ps} for details).
Unlike the former, the latter has plenty of convenient parameterizations.
The only disadvantage of this description is that the
space  of  connections  is  infinite  dimensional.  In  this  paper  (sect.
\ref{fdm}) we
consider an alternative description of the moduli space in which the  role
of the space of all connections on a Riemann surface is played by a  finite
dimensional manifold. The idea is quite familiar both  from  lattice  gauge
theory and from \v Cech cohomology. Namely, consider a  triangulation  of  a
compact Riemann surface $S$ (with boundary, in  general).  Then  we  get  a
graph formed by the vertices and the edges of this triangulation. By a graph
connection (or lattice gauge field)
we mean an assignment of a group element of the gauge group $G$ to
each (oriented) edge. The group of lattice gauge transformations ${\cal G}^l$
acting on the space of graph connections in a natural way is simply a
product of several copies of $G$, one copy for each vertex of  the  graph.
A  flat graph connection  satisfies the condition that  the
monodromies around all the faces of the triangulation are equal to $1\!\in\!G$.
(The monodromy is the product of group elements corresponding to  the
consecutive edges of a face, whatever shape
of faces is used.  One has to take account only of the orientations of the
edges in an obvious way.) Now, it is a standard assertion that the moduli
space of (smooth)  flat connections on $S$ is isomorphic to the space of flat
graph connections modulo graph gauge transformations. (This is in fact
nothing but the statement in  \v Cech cohomology that this  space is
represented by $H^1(S,G)$.) Dealing with a surface with holes amounts to
saying that some faces of the triangulation are left  empty and one does not
have to require anything about the corresponding monodromies. It is important
to note that if a graph $l$ is obtained from a triangulation of a surface $S$
it can be endowed with  an  additional structure which (together with the
graph itself) contains all the information  about the topology of the
surface. We suppose that $S$ is oriented. The orientation of $S$ induces a
cyclic order of the ends of edges incident to each vertex. A graph $l$ with a
given cyclic order at each vertex is called a  fat  graph.  If $S$ has at
least one hole the most economical way is to consider  a  fat graph with all
the faces empty, what is always possible. Conversely,  given a fat graph $l$
the corresponding surface can be restored by replacing edges of $l$ by strips
glued together at vertices respecting the cyclic order (cf. \mbox{fig. 1}).
Summarizing, in order to describe the  moduli  space ${\cal M}$ of flat
connections on a surface $S$ with holes we choose a fat graph corresponding
to $S$ (this choice is not unique) and consider the quotient of the space of
graph connections ${\cal A}^l$ by  the  action  of  graph  gauge
transformations,  ${\cal  M} = {\cal A}^l/{\cal G}^l$.

Having described the moduli space as a manifold we are  interested  now  in
describing its Poisson structure. Let us forget for a moment  that  we  can
define a Poisson structure on ${\cal M}$ by reduction of the space of  all
(smooth) connections on $S$ and try instead to define a  Poisson  structure
on ${\cal A}^l$ in such a way that it can be pulled down on ${\cal M}$
~\footnote{As it will be proved in sect. \ref{fdm} we obtain  in  this  way
the same Poisson structure as defined by the reduction procedure
from smooth connections.}.
We would like to have such a Poisson
structure on ${\cal A}^l$ that the projection ${\cal A}^l \rightarrow
{\cal M}$ will be a Poisson map. This can be achieved if ${\cal G}^l$ will act on
${\cal A}^l$ in a Poisson way (see ref. \cite{STS} for the  definition  of
Poisson group actions
on Poisson manifolds). For this aim we have to define  first
a Poisson-Lie structure on ${\cal G}^l$ itself. The group of graph gauge
transformations ${\cal G}^l$ is the direct product of several copies of $G$,
with one copy per each vertex of $l$.
Let us define the Poisson structure on ${\cal G}^l$ as a direct
product of Poisson structures on each copy of $G$ in ${\cal G}^l$. The
latter can be defined  independently at each vertex. (To define a Poisson
structure on $G$  one  has to choose a classical $r$-matrix.) Now we look for
a Poisson structure on ${\cal A}^l$. The requirement that the action of
${\cal G}^l$
is a Poisson one  is almost sufficient to determine the Poisson structure on
${\cal  A}^l$.  The
ambiguity amounts in fact to choosing a linear order of  ends  of  edges at
each vertex. Therefore, instead of fat graphs we have to  deal  with  graphs
with linear order. Let us call such graphs {\em ciliated fat graphs}.
A ciliated fat graph can be considered as a fat graph  with  an  additional
structure
(the fat graph underlying a given ciliated fat one is  restored  uniquely).
This additional structure (linear order at each vertex) can be  represented
by picturing the underlying fat graph on a sheet of paper  in  such  a  way
that the cyclic order is everywhere, say, counterclockwise and by placing
at each vertex a small cilium separating the minimal and the maximal end
incident to that vertex. As it was mentioned, a fat graph defines a
surface, that is an oriented surface with holes \mbox{(fig. 1)}; a
ciliated fat graph, similarly, defines a ciliated surface, that is
an oriented surface with holes and
with some points marked on the boundary (\mbox{fig. 2}). Thus for every ciliated
fat graph we  have an associated Poisson manifold, namely the space of
graph connection  endowed with  an $r$-matrix Poisson structure. It may happen
of course that two different ciliated graphs give isomorphic Poisson
manifolds of  graph  connections. In particular, one can show that the
isomorphism class of the arising Poisson manifold depends only on the
diffeomorphism class of the corresponding  ciliated surface.

     It may be worth mentioning some distinguished examples of graphs and
corresponding Poisson manifolds. The Poisson manifold corresponding to the
graph consisting only of two vertices and one edge (\mbox{fig. 3a}) coincides
with the Poisson-Lie group $G$ provided the $r$-matrices chosen at the vertices
are related by the operation of permutation of tensor factors
($r=r_{12} \mapsto r_{21}$).
 With the same condition on $r$-matrices, the
graph consisting of two vertices and two edges connecting  them  (\mbox{fig. 3b})
yields the manifold $G \times G$ endowed with a Poisson-Lie structure
coinciding with that of the double $D \simeq G \times G$. If we take the
same $r$-matrices at two vertices we get $D_+$ as our Poisson manifold
(see ref. \cite{STS} for the definitions of doubles). Finally, the graph
consisting of
one vertex and one edge (\mbox{fig. 3c}) corresponds to the Poisson manifold
   $G^*$, the dual Poisson-Lie group.

     The following operations with graphs are important to discuss: $i$)
erasing an edge (\mbox{fig. 4}), $ii$) contracting an edge (\mbox{fig. 5}),
$iii$) gluing graph(s) (\mbox{fig. 6}),
and $iv$) adding a loop (see sect. \ref{fdm}).
The linear orders at the vertices
touched by such an operation descend from those of the original graph
in a more or less obvious way (cf. \mbox{figs. 4,5,6}). We have to mention
only that there are in fact two ways to contract an edge which differ in what
happens to the cilia. The operation of gluing deserves some explanation.
Given two vertices on a graph with the same number $N$ of ends of edges
incident to them we can form a new graph by erasing both vertices and
gluing together thus liberated edges. (The $k$-th end liberated at one
vertex is to be glued to the $(N-k)$-th end at the other vertex.)
Note that with help of this operation one can glue together two
different graphs obtaining a single new one.

     For the operations on graphs just described there exist natural
maps between the corresponding spaces of graph connections.
These maps are in fact
projections in directions shown by the arrows in \mbox{figs. 4,5,6}.
A pleasant feature is that these maps turn out to be Poisson maps. More
precisely,  in  case  of  gluing  one  has  to  require   that   the
$r$-matrices
at two vertices to be glued are related by permutation of tensor
factors. Consider, for instance, a map corresponding to gluing together  two
simplest graphs (\mbox{fig. 7a}) each of which represents  the  Poisson-Lie
group
$G$ (an edge with two vertices). The result of gluing is again a graph
of the same shape while the corresponding map of graph connections,
$G \times G \rightarrow G$, is simply the group product which is known to be
a Poisson map. Similarly, gluing together the graphs representing $D$
gives the Poisson map $D \times D \rightarrow D$ (\mbox{fig. 7b})
corresponding to the group multiplication. Contracting one of two edges  of
the $D$ graph (\mbox{fig. 7c}) one obtains the Poisson map $D \rightarrow
G^*$.  As a Poisson manifold, the dual group $G^*$ can be identified with the
coset $D/G_\Delta$ where $G_\Delta$ is the diagonal subgroup in $D \simeq G
\times G$ (cf. ref. \cite{STS}). The isomorphism of $G^*$ with the coset
$D/G_\Delta$ shows that there is a Poisson action of $D$ on $G^*$, i.e. a
Poisson map $D \times G^* \rightarrow G^*$ which again can be described by
gluing graphs (as shown in \mbox{fig. 7d}). Looking at the pictures above
suggests the following generalisation of the notion of a double. Namely, we
can define a Poisson-Lie group, called in general a {\em
polyuble}~\footnote{We dedicate the Poisson-Lie groups of this type to
I.V.Polyubin}, by the ciliated fat graph consisting of two vertices and
several edges connecting them (analogously to the case of the double, the
$r$-matrices at two vertices should be related by the operation of
permutation of tensor factors, while the order of ends should be opposite;
\mbox{fig. 7e}). An immediate observation is
that on the space of graph connections ${\cal A}^l$ for an arbitrary
ciliated fat graph $l$ there is a Poisson action of a polyuble, $P_n(G)$,
adjusted  to  each vertex, where $n$ is the number of legs at that
vertex (see, \mbox{fig. 7f}).  Thus the space ${\cal A}^l$ is
a homogeneous space
for the group $P^l$ which is a direct product (in the sense of Poisson
groups) of $P_n(G)$'s.  Note also that the group of  graph gauge
transformations ${\cal G}^l$ which gives us the moduli space ${\cal M} = {\cal
A}^l/{\cal G}^l$ is a Poisson subgroup in $P^l$.  (Any individual polyuble,
disregarding for the moment the Poisson structure, is a product $G \times
\ldots \times G$ and contains the diagonal subgroup which turns out to be a
Poisson subgroup.)

Finally, it is worth mentioning that some particular cases of Poisson
manifolds
defined by graphs have been considered in literature. Namely the
Poisson  manifold of graph connections on a graph corresponding to the
boundary  of a polygon was suggested in ref. \cite{STS} as a discrete
approximation of current algebra coadjoint space.
(See also refs. \cite{AFSV,AFS} where this discrete approximation was
used to investigate WZNW conformal model.)

\section{Poisson structure of moduli spaces} \label{ps}

In this section we shall describe a Poisson structure on the space of
flat connections modulo gauge transformations on Riemann surfaces with
holes by means of a reduction of the space of all smooth connections
on them.

Let $S$ be an oriented compact Riemann surface with holes. Let ${\cal A}$ be
the space of smooth $G$-connections on it, where $G$ is a
reductive complex Lie
group with the Lie algebra ${\frak g}$ with a chosen nondegenerate
invariant  quadratic  form which we denote by ${\rm tr}$.
The space ${\cal A}$ is in a natural way a
symplectic manifold with the symplectic structure
\begin{equation}\label{Omega}
    \Omega = \int_{S} {\rm tr}\, (\delta A \wedge \delta A),   \label{Symp}
\end{equation}
where $A \in {\cal A}$ is a ${{\frak g}}$-valued 1-form on $S$, $\delta$
is the external differential on ${\cal A}$,
and $\wedge$ is a shorthand way to denote the wedge product both on
${\cal A}$ and on $S$.
   This symplectic structure is well known to be invariant with
respect to the gauge transformations
\begin{equation}
   A \mapsto g^{-1}A\,g + g^{-1}dg,
\end{equation}
where $g$ is a $G$-valued function on $S$.

Now let us try to define the momentum mapping  for  this  action.  One  can
easily check that the infinitesimal gauge  transformation  $\epsilon$  is
generated by the Hamiltonian function
\begin{equation}
   H_{\epsilon} = \int_{S} {\rm tr}\, (\epsilon (dA + A \wedge A)) +
                  \int_{\partial S} {\rm tr}\, (\epsilon A). \label{Ham}
\end{equation}

The Hamiltonian generating a given transformation is defined only up  to an
additive constant and therefore the Poisson brackets between them, in
general, reproduce the commutation relations between the elements of the
gauge algebra only up to a cocycle:  \begin{equation} \{
   H_{\epsilon_1},H_{\epsilon_2} \} = H_{[\epsilon_1,\epsilon_2]} +
   c(\epsilon_1,\epsilon_2).
\end{equation}

In our case
\begin{equation}
c(\epsilon_1,\epsilon_2) = \int_{\partial S} {\rm tr}\, (\epsilon_1 d\epsilon_2).
\label{MC}
\end{equation}

One can prove that this cocycle is nontrivial and therefore we  can  define
the momentum mapping not for the algebra of gauge transformations itself,
but  only  for  its
central extension by the 2-cocycle eq.(\ref{MC}).

Let ${\frak g}^S$ denote the algebra of gauge transformations centrally
extended by (\ref{MC}) and let ${\cal G}^S$ be the corresponding group.
The space ${\frak g}^S$ is the
space of pairs $(\epsilon, z)$, where $\epsilon$ is a ${\frak g}$-valued
function and $z$ is a complex number. Let us
consider the space $({\frak g}^S)^{*}$ consisting of triples $(R,B,x)$ where
$R$  is a ${\frak g}$-valued two form on $S$, $B$ is a ${\frak g}$-valued
1-form on the  boundary of $S$ and $x$ is a complex number.  There is a
nondegenerate pairing $<\,,>$ between ${\frak g}^S$ and $({\frak g}^S)^{*}$,
\begin{equation}
   <(R,B,x),(\epsilon,z) > \, = \int_{S} {\rm tr}\, (\epsilon R) +
                       \int_{\partial S} {\rm tr}\, (\epsilon B)
                      + zx.
\end{equation}

The momentum map for the action of
${\frak g}^S$ can be defined now as a mapping ${\cal A} \rightarrow
({\frak g}^S)^{*}$,
given by the
curvature and by the restriction of the connection form to the boundary.
\begin{equation}
 A \mapsto ( dA + A \wedge A, \,A\vert_{\partial S} , \,1 )
\end{equation}

Now consider the Hamiltonian reduction of ${\cal A}$ with  respect  to
${\cal G}^S_0$, the group of gauge transformations equal to the identity on
the boundary which yields the space of flat connections on $S$ modulo gauge
transformations from ${\cal G}^S_0$,
\begin{equation}
{\cal M}_0 = \{A \in {\cal A} \,\vert\, dA+A\wedge A=0\}\,/\,{\cal G}^S_0.
\end{equation}
The space ${\cal M}_0$ can also be considered as the space of values of 
flat connections restricted to the boundary. It is well known that  the
space of $G$-connections on a  circle  can  be  identified  with  the
coadjoint space of the affine Kac-Moody algebra with the standard
Kirillov-Kostant Poisson structure. The following proposition shows that
these two Poisson structures are related:
\begin{proposition} The mapping
from the space ${\cal M}_0$ to the Kac-Moody  coadjoint representation space
sending a flat connection on the Riemann surface $S$ to its restriction to a
component of the boundary is a Poisson mapping.
\end{proposition}
{\em Proof. } This mapping is essentially the momentum mapping for the
action of gauge transformations. $\Box$

Now let us consider the quotient of the space ${\cal M}_0$ by the whole
group ${\cal G}^S$ (the group ${\cal G}^S$ acts on ${\cal M}_0$  because the
group ${\cal G}^S_0$ of gauge transformations equal to the identity
on  the  boundary  is normal in ${\cal G}^S$). The quotient space,
\begin{equation}
{\cal M} = \{A \in {\cal A} \,\vert\, dA+A\wedge A=0\}\,/\,{\cal G}^S_1,
\end{equation}
is a finite dimensional Poisson manifold. Its symplectic leaves are in
one-to-one correspondence with the coadjoint orbits of the centrally
extended  group of gauge transformations  which
in turn are  parameterized by the conjugacy classes of monodromies
around the holes.  Thus we have
\begin{proposition}\label{prop2}
The space of flat $G$-connections modulo gauge transformations, ${\cal
M}$, on a Riemann surface with holes inherits a Poisson structure from the
space of all (smooth) $G$-connections. The symplectic leaves of  this  structure  
are
parameterised by the conjugacy classes of monodromies around holes.
\end{proposition}

\section{Graph connections} \label{fdm}

In this section we shall construct a Poisson structure on the space of
graph connections, ${\cal A}^l$, in such a way that the lattice  gauge
group, ${\cal G}^l$, endowed with a nontrivial $r$-matrix Lie-Poisson  structure
   acts  on ${\cal A}^l$ in a Poisson way.

Let $l$ be a ciliated
fat graph homotopically equivalent to a Riemann surface $S$
with holes. Denote by $E(l)$ the set of ends of edges of $l$ and by  $N(l)$
the set of its vertices. Each element of $N(l)$ corresponds to
the subset of $E(l)$ of ends of edges incident to a given vertex. In what
follows we shall identify elements of $N(l)$ with  the  corresponding
subsets. A mapping which sends an end of an edge $\alpha$ to the opposite
end of the same edge $\bird{\alpha}$ is an involution of the set $E(l)$.
The ciliated fatness of $l$ defines an ordering inside
each $n \in N(l)$.
One can easily see that such data -- a set divided into ordered
subsets and an involution of it without fixed points -- unambiguously
define a ciliated fat graph.
Let $[\alpha]$ be the vertex containing
$\alpha$ and $[\alpha,\bird{\alpha}]$ be the edge linking
$\alpha$ and $\bird{\alpha}$.

   Let us call a {\em graph connection} on a graph $l$
an assignment of an element
${\bf A}_{\alpha}$ of the group $G$ to each $\alpha \in E(l)$ such that
\footnote {Perhaps it would be more natural to assign a group
element to each edge, as we did in sect. \ref{cfg}  above, rather than to
each end of an edge. However, in this case we would have to choose some
orientations of the edges. Then we would have to have a definition of the
Poisson manifold ${\cal A}^l$ which would depend on an {\em oriented}
ciliated fat graph. In such a case it would be possible to prove that two
Poisson manifolds corresponding to two graphs differing only by
their orientations are isomorphic. We prefer to get rid of this
complication at the price of a slightly more complicated notation.}

\begin{equation} \label{Inv}
  {\bf A}_{\alpha\!^\vee} = {\bf A}_{\alpha}^{-1} .
\end{equation}

The lattice gauge group  ${\cal G}^l$  is  a
product of finite dimensional groups $G$ --- one  copy for  each  vertex  of
the graph. The group ${\cal G}^l$ acts on ${\cal A}^l$ in a natural way:
\begin{equation}
{\bf A}_{\alpha} \mapsto {\bf g}_{\alpha\!^\vee}^{-1}\,
{\bf A}_{\alpha} \,{\bf g}_\alpha .
\end{equation}

 The space of graph connections can be considered as a  quotient  space  of
the space of flat connections on a surface $S$. Indeed, let us take the
surface $S$ corresponding to the graph $l$ and imbed the graph into it in
a way such that $S$ is contractable to the image of $l$. Then for a
(smooth) connection $A$ on $S$ we can construct a graph connection on $l$
assigning to $\alpha \in E(l)$ the parallel transport operator along the edge
linking $\bird{\alpha}$ and $\alpha$. This graph connection does not change
if we transform the connection $A$ by a gauge transformation equal to the
identity at the vertices. It is clear that every graph connection can be
continuously extended to the surface  and therefore the space of graph
connections ${\cal A}^l$ can be represented as a quotient,

\begin{equation}
{\cal A}^l \cong \{A \in {\cal A} \,\vert\, dA+A\wedge A=0\}\,/\,{\cal G}^S_1,
\end{equation}
where ${\cal G}^S_1$ is the
group of gauge transformations equal to the identity at the
vertices. Of course, this representation is defined only up to
the action of the graph gauge group and, therefore, the isomorphism
between the spaces ${\cal M}$ and ${\cal A}^l/{\cal G}^l$ is canonical.

  This isomorphism shows that although
the space ${\cal A}^l$ has so far no {\em a priori} Poisson structure,
the space ${\cal A}^l/{\cal G}^l$ has one. Our aim is to introduce a Poisson
structure on ${\cal A}^l$ compatible with that on ${\cal A}^l/{\cal G}^l$ and
with the graph gauge group action.

Let us fix for each vertex $n$ of the graph
a classical $r$-matrix $r(n) \in {\frak g} \otimes {\frak g}$
, that is to say, a solution
of the classical Yang-Baxter equation:
\begin{equation} \label{YB}
[r_{12}(n),r_{13}(n)] + [r_{12}(n),r_{23}(n)] + [r_{13}(n),r_{23}(n)] = 0
\end{equation}
such that
\begin{equation}\label{symmpart}
          \frac{1}{2}( r_{12}(n) + r_{21}(n)) = t,
\end{equation}
where $t \in {\frak g} \otimes {\frak g}$ is a quadratic Casimir element:
\begin{equation}
t = \sum e_i \otimes e_i,
\end{equation}
where $\{e_i\}$ is an orthonormal basis in ${\frak g}$
~\footnote{Note that although the $r$-matrix, $r(n)$, is allowed
to differ for different vertices, its symmetric part, $t$, is required to be
the same everywhere. }.

Let us define a bivector field $B$
on ${\cal A}^l$ as
\begin{equation} \label{Pois}
B = \sum_{n \in N(l)}\left(\sum_{\alpha ,\beta \in n; \alpha < \beta}
          r^{ij}(n)X_i^{\alpha} \wedge X_j^{\beta} +
 \frac{1}{2}\sum_{\alpha \in n} r^{ij}(n)X_i^{\alpha}
\wedge X_j^{\alpha}\right),
\end{equation}
where $X_i^{\alpha} = L_i^\alpha - R_i^{\alpha\!^\vee}$, $L_i^\alpha$ and
$R_i^\alpha$ are, respectively, the left- and right-invariant vector fields
corresponding to the element $e_i \in {\frak g}$
on the group assigned to $\alpha \in E(l)$
and $r^{ij}(n)$ is the $r$-matrix at the vertex $n$ written in
the basis $\{e_i\}$. Note that the vector  fields  $X_i^{\alpha}$
are chosen to be consistent with eq.(\ref{Inv}).
\begin{proposition}

a) The bivector $B$ defines a Poisson structure on ${\cal A}^l$.
b) The group ${\cal G}^l$ endowed with the direct product
Poisson-Lie structure acts on ${\cal A}^l$ in a Poisson way.
\end{proposition}

The proof can be obtained by a straightforward check.

\vspace{4mm}

Sometimes it is however more convenient to use other ways of presenting
the Poisson bivector (\ref{Pois}). If one separates explicitly the
symmetric, $t$, and skew-symmetric, $r_a = \frac{1}{2}(r_{12}-r_{21})$,
parts of the r-matrix, so that $r=r_a+t$, one gets
\begin{equation} \label{Pois1}
B = \sum_n \left( r_a^{ij}(n)X_i^{\Delta}(n)\otimes X_j^{\Delta}(n)  +
    \sum_{\alpha,\beta \in n}(n,\alpha,\beta)
     \sum_iX_i^{\alpha}\otimes X_i^{\beta}\right),
\end{equation}
where $X_i^{\Delta}(n) = \sum_{\alpha \in n}X_i^{\alpha}$ and
\begin{equation}
(n,\alpha,\beta) =
\left\{ \begin{array}{ll} 1 & \alpha > \beta \\
                          0 & \alpha = \beta \\
                         -1 & \alpha < \beta
\end{array}\right.,\qquad \mbox{ for } \alpha,\beta \in n.
\end{equation}

Since the vectors $X_i^{\Delta}(n)$ are tangent to ${\cal G}^l$-orbits,
one sees that the Poisson bracket induced by eq.(\ref{Pois}) on the quotient
${\cal M} = {\cal A}^l/{\cal G}^l$ does not change if the skew-symmetric part
$r_a$ of the $r$-matrix is changed.

Another way of defining the Poisson structure is to give explicit expressions
for the Poisson brackets between matrix elements
of ${\bf A}_{\alpha}$ in some representation of the group $G$;
we consider these matrix elements as functions on ${\cal A}^l$.
(We shall denote matrices representing ${\bf A}_{\alpha}$ and $r$
by the same symbols, ${\bf A}_{\alpha}$ and $r$, respectively.)
\begin{equation}\label{edge}
\{{\bf A}_{\alpha} \otimes {\bf A}_{\alpha} \} =
r_a(1)\,({\bf A}_{\alpha} \otimes {\bf A}_{\alpha})  +
({\bf A}_{\alpha} \otimes {\bf A}_{\alpha})\,r_a(2)
\end{equation}
for the case $[\alpha] \neq [\alpha^\vee]$. Here $r(1) = r([\alpha])$,
$r(2) = r([\alpha^\vee])$.
\begin{equation}\label{loop}
\{{\bf A}_\alpha,{\bf A}_\alpha\} = r_a\,({\bf A}_\alpha\otimes {\bf
A}_\alpha) + ({\bf A}_\alpha\otimes {\bf A}_\alpha)\,r_a + (1 \otimes {\bf
A}_\alpha)\,r _{21}\,({\bf A}_\alpha \otimes 1) -
({\bf A}_\alpha \otimes 1)\,r\,(1 \otimes {\bf A}_\alpha)
\end{equation}
for the case $[\alpha] = [{\alpha}^\vee]$,\ $\alpha < \alpha^\vee$.
\begin{equation}\label{twoedges} \{{\bf A}_{\alpha} \otimes {\bf A}_{\beta}
\} = r\,({\bf A}_{\alpha} \otimes {\bf A}_{\beta})
\end{equation}
for the case
$[\alpha] = [\beta] \neq [\alpha^\vee] \neq [\beta^\vee] \neq [\alpha]$,
$\alpha < \beta$.
\begin{equation}\label{double}
\{{\bf A}_\alpha,{\bf A}_\beta\} =
r(1)\,({\bf A}_\alpha \otimes {\bf A}_\beta) +
({\bf A}_\alpha \otimes {\bf A}_\beta)\,r(2)
\end{equation}
for the case
$[\alpha] = [\beta] \neq [\alpha^\vee] = [\beta^\vee]$,\ $\alpha < \beta$,\
$\alpha^\vee < \beta^\vee$,\ $r(1) = r([\alpha])$, $r(2) = r([\alpha^\vee])$.
\begin{equation}\label{torus}
\{{\bf A}_\alpha,{\bf A}_\beta\} = r\,({\bf A}_\alpha\otimes {\bf A}_\beta) +
({\bf A}_\alpha\otimes {\bf A}_\beta)\,r +
(1 \otimes {\bf A}_\beta)\,r_{21}\,({\bf A}_\alpha \otimes 1)-
({\bf A}_\alpha \otimes 1)\,r\,(1 \otimes {\bf A}_\beta).
\end{equation}
for the case
$[\alpha] = [\beta] = [\alpha^\vee] = [\beta^\vee]$ and $\alpha <
\beta < \alpha^\vee < \beta^\vee$.  Unfortunately the complete list of all
possible configurations of one or two edges and cilia is rather long (there
are fourteen of them) and we stop here. The reader can easily observe how one
can write down expressions for other configurations by analogy.

\vspace{4mm}

As it was described in sect. \ref{cfg}, there exist such operations on
ciliated fat graphs as erasing an edge, contracting an edge towards a vertex,
gluing two vertices of the same valence and adding a loop. One can also
change a graph to another one corresponding to the same ciliated surface. All
these transformations induce mappings between the corresponding spaces of
graph connections.  Let us now describe them explicitly.

{\em Erasing an edge (\mbox{fig. 4}).} This operation is the most obvious one.
The mapping between graph connections is just given by forgetting the group
element assigned to the edge to be erased.

{\em Contracting an edge (\mbox{fig. 5}).} This operation can be applied to
an edge
with distinct ends (i.e. $[\alpha] \neq [\alpha^\vee]$). Let $\alpha$ be an
end of such an edge.  (In \mbox{fig. 5}, it is the right one for the
projection $R$ and the left one for $L$.) Make ${\bf A}_\alpha$ be equal to
identity by applying a gauge transformation (that is the action of one copy
of $G$) at the vertex $[\alpha^\vee]$.  Erase the cilium at the vertex
$[\alpha^\vee]$.  Then contract the edge $[\alpha,\alpha^\vee]$ leaving the
group elements on the other edges unchanged (as they were after the above
gauge transformation).

Note that, as it is shown in \mbox{fig. 5}, this operation depends not only
on the edge but also on the choice of a particular end of it. To emphasize
this we say that we contract the edge towards a vertex, in our case
$[\alpha]$.

{\em Gluing two vertices (\mbox{fig. 6})}. This operation can be applied to two
vertices $n$ and $n'$ having the same valence, i.e. $|n|=|n'|$, and such
that their $r_a$-matrices are opposite, i.e. $r_a(n)=-r_a(n')$.
Disconnect
the ends of edges at the vertices and connect them in the order prescribed by
gluing (\mbox{fig.6)} inserting an arbitrarily ciliated 2-valent vertex at
each connection. Until now we left the group elements on the edges unchanged.
Now take each inserted 2-vertex and contract towards it one of the two
incident edges.

{\em Adding a loop.} One can add a loop (an edge $[\alpha, \alpha^\vee]$
with $[\alpha = \alpha^\vee]$ to a vertex between two consecutive ends of
edges.  Assign the unit group element to the new loop.

{\em Ciliated graphs and ciliated surfaces.}
As it was mentioned several times above, a graph imbedded into an oriented
surface inherits a fatness (cyclic order of ends of edges at vertices).
Assume now that we have a graph imbedded into a surface in such a way that the
vertices are mapped into the boundary. This graph inherits a ciliated
fatness since there is a canonical linear order of the ends of edges meeting
at a boundary point. On the other hand, given a ciliated fat graph imbedded
into the corresponding surface (that is, we assume that the surface is
retractable to the image of the graph) there exists a unique up to the
isotopy way to move its vertices to the boundary reproducing the given
ciliation. We have just to move each vertex to the boundary component
which the cilium looks onto. If we now erase the edges of the graph and leave
the cilia stuck out off the boundary components  we get
a {\em ciliated surface} (e.g., \mbox{figs. 2b, 2d}).

Suppose now we have two ciliated fat graphs $l$ and $l'$ corresponding
to the same ciliated surface. (This means, in particular, that their vertices
are identified.)
We are going to construct an isomorphism
${\cal A}^l \stackrel{\sim}{\longrightarrow} {\cal A}^{l'}$
between the spaces of graph connections on them.
Let $\alpha \in E(l')$ be an end of an edge of $l'$. Take the edge
$[\alpha, \alpha^\vee]$ of $l$ and retract it to the graph $l'$. We obtain
a path on the graph $l$ connecting the vertices $[\alpha]$ and $[\alpha^\vee]$
and isotopic to the edge $[\alpha, \alpha^\vee]$. Assign to $\alpha$ the monodromy
of the graph connection ${\cal A}^l$ along this path. Carrying out this
procedure for all $\alpha \in E(l')$ we get the desired isomorphism.
\vspace{4mm}

Now let us summarize some properties of the spaces of graph connections
equipped with the Poisson bracket eq.(\ref{Pois}).
\begin{proposition}\label{opers}~\newline
1) The mappings between graph connections corresponding to erasing an edge,
contracting an edge towards a vertex and gluing two vertices are the
Poisson projections onto the image.
\newline
2) The isomorphism of the spaces of graph connections for two graphs
corresponding to isomorphic ciliated surfaces is an isomorphism of Poisson
manifolds.
\newline
3) Adding a loop is a Poisson imbedding.
\end{proposition}

The proof of the proposition is a straightforward and not very complicated
explicit calculation that we omit here. Let us mention only the following
statement useful for the proof as well as by itself.

Let $f$ be a face of a ciliated fat graph $l$ such that there are no cilia
looking into $f$.  Let ${\cal A}^l(h,f)$ be the set of graph connections with
the monodromy around the face $f$ of $l$ conjugated to $h \in G$. Then
${\cal A}^l(h,f)$ is a Poisson submanifold in ${\cal A}^l$.

\vspace{5mm}
Let us proceed now to the relation between the space of graph
connections and the space of ordinary connections.
\begin{proposition}
The quotient of the space of graph connections by the graph gauge group
is isomorphic as a Poisson manifold to the quotient of  the  space  of
flat connections on the corresponding Riemann surface by the gauge group,
i.e.,
\begin{equation}\label{Isom}
{\cal A}^l/{\cal G}^l \cong {\cal M}.
\end{equation}
\end{proposition}

{\em Remark 1.}
Let us note that this statement shows that all
the ambiguities in the construction of the space ${\cal A}^l$ -- such  as
choices of ordering and of $r$-matrices -- do not influence the Poisson
structure of its quotient by the gauge group.
The latter depends only on the cyclic order and on the symmetric part, $t$, of
the $r$-matrices (cf. eq.(\ref{symmpart})). This could not be otherwise,
because these are just the data defining the Poisson manifold ${\cal M}$ by
eq.(\ref{Omega}), provided the surface $S$ there is defined by the ciliated fat
graph $l$ here and the invariant scalar product ${\rm tr}$ there is defined by
the Casimir element $t$ here.
However it is impossible to introduce  a Poisson
structure on ${\cal A}^l$ compatible with that on the gauge  quotient without
fixing nontrivial $r$-matrices. Note also that topologically these moduli
spaces are always isomorphic to products of several copies of the group $G$
modulo the overall $G$-conjugation, although they are  not  isomorphic to each
other as Poisson manifolds. For example a sphere with three holes and a torus
with one hole give topologically the same spaces, $(G \times G)/ {\sl Ad}\,G$,
while
the Poisson structure for, e.g., $G = SL(2)$ is trivial in the first case and
nontrivial in the second one.

{\em Remark 2.} The description of the moduli space ${\cal M}$ of flat connections 
in the graph language has an advantage that this language allows us to describe 
rather explicitly the space of functions on ${\cal M}$ using representation theory. 
In particular one can construct a linear basis in the space of regular functions on 
${\cal M}$ in the following way.

Assign an irreducible representation
$\pi_{\alpha}$ of $G$ in a space $V_\alpha$ to each $\alpha \in E(l)$ in
such a way that
$\pi_{\bird{\alpha}}  =  \pi_{\alpha}^*$
and assign an intertwiner $C_n \in \mbox{Inv}
(\otimes_{\alpha \in n}V_{\alpha}^*)$ to each vertex $n$.
We can consider matrices from $\mbox{End\,}V_{\alpha}$ as belonging to
$V = \otimes_{\alpha \in E(l)}V_\alpha$ and the
intertwiners $C_n$ as belonging to its dual, $V^*$.

For each such data  $(l,C_\bullet,\pi_\bullet)$  we can define
a function $\psi(l,C_\bullet,\pi_{\bullet})$ on ${\cal A}^l$
\begin{equation}
\psi(l, C_\bullet,\pi_\bullet)(\{{\bf A}_{\alpha}\}) =
 <\bigotimes_n C_n, \bigotimes_{\alpha \in E_1(l)}\pi_{\alpha}({\bf A}_{\alpha})>
\end{equation}
where $E_1(l) \subset E(l)$ is a set of ends of  edges  containing  exactly
one end of each edge. The ambiguity in  the choice of this  set  is inessential
because $\pi_{\bird{\alpha\!}}({\bf A}_{\bird{\alpha\!}}) =
\pi_\alpha({\bf A}_{\alpha})$ as an element of $V_\alpha \otimes V_\alpha^*$.

One can easily verify that all such functions are ${\cal G}^l$ invariant and that 
they indeed form a complete set of functions on ${\cal M}$. The latter is an 
obvious consequence of the Peter-Weyl theorem. 

\vspace{5mm}
{\em Proof of the proposition 5.} 
The Poisson bracket of two functions $\Psi$ and $\Phi$ on the space 
${\cal A}$ of smooth connections on $S$ can be written as
\begin{equation}
\{\Psi,\Phi\}_{_{\scriptstyle S}} = \int_S {\rm tr}\, \left( \frac{\delta 
\Psi}{\delta A}
            \wedge \frac{\delta \Phi}{\delta A}\right),\label{smooth}
\end{equation}

To prove the proposition we need to compute the Poisson bivector on ${\cal M}$ 
induced from eq.(\ref{smooth}) by the reduction procedure described in the sect. 3 
and compare the result with the bivector induced by eq.(\ref{Pois}).

In order to be able to work with the bracket eq.(\ref{smooth}) and build a bridge 
between the smooth and the combinatorial approaches to the Poisson brackets on flat 
connections let us first compute the Poisson bracket using the formula 
(\ref{smooth}) in one particular case. Let $I_1$ and $I_2$ be two oriented 
intervals imbedded into $S$ and intersecting transversally. Let us compute using 
eq.(\ref{smooth}) the Poisson bracket between two arbitrary functions $\Psi$ and 
$\Phi$ of the corresponding monodromies considered as functions on ${\cal A}$.

The result of the computation is a function on the space ${\cal A}$. However, for 
our further purposes, we need to compute only the restriction of the result to the  
connections such that their restrictions to the intervals vanishes everywhere 
except for two subintervals containing the ends of $I_1$ and $I_2$ and none of 
their intersection points. In this case the expression for the bracket is 
especially simple, it depends only on the monodromies along the segments and can be 
straightforwardly computed from eq.(\ref{smooth}):
\begin{equation}
\{\Psi,\Phi\} = t^{ij} (R_i\Psi) (R_j\Phi)\sum_{k \in I_1 \cap I_2}\varepsilon(k) . 
\label{P1}
\end{equation}
Here $k$ runs over the intersection points, $\varepsilon(k)$ is $1$ or $-1$ 
if the first segment crosses the second one from the left or from the right, 
respectively, $\{R_i\}$ is a basis of the left-invariant vector fields on $G$ and 
$t^{ij}$ is the matrix of the quadratic Casimir $t \in {\frak g} \otimes {\frak g}$ 
(e.g., eq. (\ref{symmpart})).

\vspace{3mm}
Note that this formula does not give Poisson brackets between functions of 
monodromy along a single segment. Moreover the formula (\ref{smooth}) is not 
applicable to compute such brackets.

Let us recall the definition of the Hamiltonian reduction 
in the language of Poisson brackets. Let ${\sf M}$ be a 
symplectic manifold with symplectic action of the group ${\sf G}$, $\mu$ be 
a momentum map (corresponding to the action of ${\sf G}$ or of a subgroup 
of it , ${\sf M}_0 = \mu^{-1}(0)$,\ 
${\sf N} = {\sf M}_0/{\sf G}$ be the reduced space and 
$\pi:{\sf M}_0 \rightarrow {\sf N}$ be the canonical projection. The 
Poisson bracket of two functions $\Psi$ and $\Phi$ on ${\sf N}$ is defined as 
follows. Let $\Psi^*$ and $\Phi^*$ be any two functions on ${\sf M}$ 
such that their restrictions to ${\sf M}_0$ coincide with $\pi^*\Psi$ 
and $\pi^*\Phi$ respectively.  
Then $\{\Psi,\Phi\}({\sf x}) := \{\Psi^*,\Phi^*\}({\sf y})$, 
for any ${\sf x} \in {\sf N}$ and any ${\sf y} \in \pi^{-1}({\sf x})$. (Note 
that this procedure includes at least three arbitrary choices: the choice 
of functions $\Psi^*$ and $\Phi^*$ for given $\Psi$ and $\Phi$ and the 
choice of ${\sf y}$ for given ${\sf x}$. We are going to make these choices  
in a way maximally simplifying the calculations.) 

This definition can be applied to our situation. We have the space of all 
connections, 
${\cal A}$, as ${\sf M}$, the space of flat connections as ${\sf M}_0$
and the moduli space, ${\cal M}$, as ${\sf N}$. Our task is to compute 
Poisson bracket on ${\cal M}$ or, equivalently, between 
${\cal G}^l$-invariant functions on ${\cal A}^l$ and, then, compare the result with 
the one given by eq.(\ref{Pois}). 

Let $\Psi$ and $\Phi$ be arbitrary two such functions and let ${\it \bf l}$ and 
${\it \bf l'}$ be two imbedding of the graph $l$ into the surface 
such that the images of vertices are disjoint and the images of edges 
are transversal. Using the mappings ${\cal A} \rightarrow {\cal A}^l$ 
given by monodromies along the edges  we can lift $\Psi$ and $\Phi$ using, 
respectively, ${\it \bf l}$ and ${\it \bf l'}$ to a ${\cal G}^S$-invariant 
functions $\Psi^*$ and $\Phi^*$ on ${\cal A}$.

Now to compute the bracket between $\Psi^*$ and $\Phi^*$ we need only to apply 
the lemma to all intersecting edges of ${\it \bf l}$ and ${\it \bf l'}$. 
To simplify the computations one can choose a convenient pair of graph 
imbeddings as well as a convenient flat connection within the given 
${\cal G}^l$-orbit.
(In fact, we need two imbeddings ${\it \bf l}$ and ${\it \bf l'}$ since the 
formula (\ref{smooth}) is not applicable for computing brackets between 
functions given by one imbedding.)

Fix a ciliation on $l$ and imbed $l$ in the surface in the way such 
that all vertices map to the boundary and all cilia look outside the 
surface. Thus we get our first embedding ${\it \bf l}$. To get the 
second imbedding ${\it \bf l'}$ deform the imbedding  ${\it \bf l}$ 
in order to make the edges of  ${\it \bf l}$ and  ${\it \bf l'}$ 
transversal and the formula (\ref{P1}) applicable. Fix a point at the 
middle of each edge. Then move each vertex along the boundary component 
a little to the left (if viewed from outside)  together with incident 
edges keeping the middle points stable and making the number of the intersection 
points between deformed and initial edges as low as possible. Such a 
deformation is illustrated in \mbox{fig. 8}.

We have one intersection point for any 
two ends of edges $\alpha \in E({\it \bf l})$ and $\beta \in E({\it \bf l'})$ 
belonging to the same vertex. Let us say that these intersection points are 
{\em associated} to this vertex. There is also one intersection point 
at the middle of each edge which we associate in an arbitrary way to 
one of the vertices of the corresponding edge.

Let us choose now a convenient connection within a given ${\cal G}^l$ orbit. 
One can fix a 
disjoint collection of patches around each vertex of ${\it \bf l}$ in such a 
way that each patch contains the corresponding vertex of ${\it \bf l'}$ as 
well as all the segments of edges between these vertices and the 
intersection points associated to them. Since the patches are disjoint 
and topologically trivial, one can make the connection  on them to be zero. 

 Note 
that since we have chosen the connection to be trivial around the vertices, 
we we can apply the formula (\ref{P1}). Note also that the 
intersection points at the middle of the edges give trivial contribution. 
Finally we get an expression for a bivector giving Poisson bracket of $\Psi$ 
and $\Phi$ as the following sum over all other intersection points
\begin{equation}
B = \sum_n \left(
    \sum_{i;\alpha,\beta \in n; \alpha < \beta}
     X_i^{\alpha}\wedge X_i^{\beta}\right),
\end{equation}
which coincides with eq.(\ref{Pois1}) up to terms vanishing on 
${\cal G}^l$-invariant functions. $\Box$

\vspace{4mm}
   In this section we described the Poisson structure on ${\cal A}^l$ which
gave us a description of the Poisson structure on ${\cal M}$ as well. As we
mentioned in the Proposition \ref{prop2} above, it is also possible to
characterize the symplectic leaves in ${\cal M}$. It might be, however,
useful to have an explicit description of the symplectic structure on those
leaves. For such a description we refer to the paper by A.Alekseev and
A.Malkin, ref.  \cite{AlMal2}, see also their work ref. \cite{AlMal1} where a
useful description of the symplectic structure on the symplectic leaves in
Poisson-Lie groups is given.


\section*{Appendix. Ruijsenaars equations}

\setcounter{equation}{0}
  In this Appendix we describe the geometric meaning of the trigonometric
Ruijsenaars Hamiltonian integrable system \cite{Ruij}, see eq. (A\ref{ruij})
below.  This system is a generalization of several integrable systems such as
rational and trigonometric Calogero system, rational Ruijsenaars system
and  finite Toda
chains. All those systems can be obtained from the trigonometric Ruijsenaars system
by suitable limiting procedures. Another aspect which makes this system
very interesting is its duality property, what means that coordinates and
Hamiltonians enter this system symmetrically, i.e., there exists an involution
of the phase space interchanging them. This property fails to be present in
all the above listed limiting cases but the rational Calogero one, where
this duality is well know even in the quasiclassical case. The quantum
version of the trigonometric Ruijsenaars system is the system of MacDonald
difference operators \cite{Md} and the duality between coordinates and
hamiltonians appears there in disguise of MacDonald's conjecture recently
proved by Cherednik \cite{Cher} by the methods quite different from those
described in the present paper. However we shall not discuss the quantum
aspects of this problem here.

  We show here that one can interpret the phase space of the trigonometric
Ruijsenaars system as a symplectic leaf of the lowest dimension in the
moduli space ${\cal M}$ of flat $G=SL(k)$
connections on a once holed torus \footnote{The relation between the 
Ruijsenaars system and moduli of flat connections on the torus was found by 
Gorsky and Nekrasov in ref. \cite{GN}.}. 
The commuting Hamiltonians described in \cite{Ruij} are certain
conjugation invariant functions of one of monodromies, the monodromy around
one of the cycles of the torus, while the coordinates are the eigenvalues of
the other.  This picture shows that the duality is nothing but the action of
the element of the mapping class group of the torus interchanging these two
cycles.

  As a by-product, we introduce a Poisson bracket, as well as a set of
commuting Hamiltonians, on an auxiliary space $G\times G$. The flows
generated by the Hamiltonians are particularly simple and the corresponding
Hamiltonian equations can be easily integrated. The projection of the Poisson
structure and the Hamiltonians on the quotient $G \times G/ {\sl Ad}\,G$ exists
and gives exactly the Ruijsenaars Hamiltonian system upon restriction to a
certain symplectic leaf. This procedure gives a way to solve the Ruijsenaars
equation explicitly.  Algorithmically, it is of course just the same as it was
proposed by S.N.M.Ruijsenaars and H.Shneider \cite{Ruij}; we just give a
natural geometric meaning to it.

Our aim is now to prove the above statement. For this purpose we have to do
the following.

1. Compute the canonical Poisson bracket on ${\cal M}=G \times G/{\sl Ad}\,G$
(using the technique developed in the main part of the paper).

2. Choose coordinates on ${\cal M}$ canonically conjugated with respect to
the Poisson bracket to the eigenvalues of one of the monodromy operators.

3. Compute a  certain function of the other
monodromy conjugacy class and verify that this gives exactly the trigonometric
Ruijsenaars Hamiltonian.

 To describe the symplectic structure on ${\cal M}$, choose the ciliated fat
graph $l$ consisting of two edges and one vertex with the ciliated fat graph
structure as shown in \mbox{fig. 9} corresponding to the once holed torus
 (\mbox{fig. 10}).

The space of
graph connections, ${\cal A}^l$, for such a graph is just a product of two
copies of the group $G$,
$$
{\cal A}^l = G \times G = \{( A, B)\},\noeq{A}
$$
where ${ A}$ and ${ B}$ are assigned to the edges of the graph as indicated
in \mbox{fig. 9}.

 The Poisson brackets on ${\cal A}^l$ are given by
the relations following from the definition, eq.(\ref{Pois}):
$$
\{{\bf A}\,,{\bf A}\} = r_a\,({\bf A}\otimes {\bf A}) +
({\bf A}\otimes {\bf A})\,r_a +
(1 \otimes {\bf A})\,r_{21}\,({\bf A} \otimes 1) -
({\bf A} \otimes 1)\,r\,(1 \otimes {\bf A}), \noeq{A}\label{td1}
$$
$$
\{{\bf B}\,,{\bf B}\} = r_a\,({\bf B}\otimes {\bf B}) +
({\bf B}\otimes {\bf B})\,r_a +
(1 \otimes {\bf B})\,r_{21}\,({\bf B} \otimes 1) -
({\bf B} \otimes 1)\,r\,(1 \otimes {\bf B}), \noeq{A}\label{td2}
$$
$$
\{{\bf A}\,,{\bf B}\} = r\,({\bf A}\otimes {\bf B}) +
({\bf A}\otimes {\bf B})\,r
+ (1 \otimes {\bf B})\,r_{21}\,({\bf A} \otimes 1) -
({\bf A} \otimes 1)\,r\,(1 \otimes {\bf B}),\noeq{A}\label{td3}
$$
where $r_a = \frac{1}{2}(r - r_{21})$; ${\bf A}$ and ${\bf B}$
are the matrix functions on $G\times G$
corresponding to $A$ and $B$, respectively, in the standard $k$-dimensional
representation.

Introduce the standard notation, $G^*$, for the group $G$ equipped with the
Poisson structure given by eq.(A\ref{td1}) and
corresponding to the graph consisting of just one loop.
The relation eq.(A\ref{td1}), which coincides, of course with
eq.(\ref{loop}), is called sometimes the {\em reflection equation.}

The projections $p_1$ and $p_2$ of ${\cal A}^l= G\times G$  onto the first
and the second factor, respectively, are obviously Poisson maps
$p_{1,2}: {\cal A}^l \rightarrow G^*$

Now let us restrict ourselves to the case of the standard $r$-matrix,
$$
   r = \sum_{\alpha > 0} E_{\alpha} \otimes E_{-\alpha} + \frac{1}{2}
\sum_i H_i \otimes H_i ~. \noeq{A}
$$
In this case one can easily derive from eqs.(A\ref{td1}--A\ref{td3}) the
following commutation relations
$$
 \{{\rm tr}\, {\bf A}^n, {\bf A}\} = 0, \;
 \{{\rm tr}\, {\bf B}^n, {\bf B}\} = 0, \noeq{A}\label{comm}
$$
$$
  \{{\rm tr}\, {\bf A}^n, {\bf B}\} = n({\bf A}^n)_0\,, \noeq{A}\label{bra}
$$
$$
  \{{\rm tr}\, {\bf B}^n, {\bf A}\} = n{\bf A}({\bf B}^n)_0\,,\noeq{A}
$$
where $({\bf X})_0$ denotes the traceless part of the matrix ${\bf X}$.
Therefore, the functions ${\rm tr}\, {\bf B}^n$ for $n=1, \ldots, k-1$
considered as Hamiltonians generate commuting flows on ${\cal A}^l$:
$$
   {\bf B}(t_1,\ldots,t_{k-1}) = {\bf B}(0,\ldots ,0), \noeq{A}\label{e1}
$$
$$
{\bf A}(t_1,\ldots,t_{k-1}) =
{\bf A}(0,\ldots ,0) \, e^{(t_1 {\bf B} + \cdots + t_{k-1}{\bf B}^{k-1})_0}\,.
\noeq{A}   \label{e2}
$$

As it was shown in the main part of the paper the lattice gauge group
${\cal G}^l$ acts on ${\cal A}^l$ in a Poisson way, and the quotient Poisson
manifold coincides with the moduli space ${\cal M}$ of smooth flat
connections on the Riemann surface corresponding to the fat graph $l$. In
our case, the group ${\cal G}^l$ is just $G$ itself (since the graph has
only one vertex) which acts on ${\bf A}$ and ${\bf B}$ by a simultaneous
conjugation:
$$
g:({\bf A},{\bf B}) \mapsto  (g{\bf A}g^{-1},\,g{\bf
A}g^{-1}).\noeq{A}
$$

 The functions ${\rm tr}\,{\bf A}^n$ and ${\rm tr}\,{\bf B}^n$ are invariant
under this action, and
therefore they descend to the moduli space ${\cal M}$ and generate commuting
flows there as well, the trajectories on ${\cal M}$ being just the projections
of those given in eqs.(A\ref{e1}),(A\ref{e2}).

However the moduli space ${\cal M}$ is in our case a Poisson manifold with
a degenerate Poisson bracket. The symplectic leaves in ${\cal M}$ correspond
to connections having a fixed conjugacy class of the monodromy around the
hole.  In our case, the latter is just the matrix ${\bf A}{\bf B}{\bf
A}^{-1}{\bf B}^{-1}$.

  Let us recall now that we are actually dealing with the case of $G=SL(k)$.
  Different symplectic leaves in ${\cal M}$ have different dimensions and the
lowest dimension among them, in this case, is $2(k-1)$. Those leaves
correspond to the monodromy around the hole being conjugated to a matrix
$x{\bf 1} + {\bf P}$,
where ${\bf 1}$ is the unit matrix, ${\rm rk}\,{\bf P} \leq 1$, and
$x \neq 0$ is a
number which parameterises this set of symplectic leaves of the lowest
dimension.  (Indeed, the only conjugation invariant of an operator of rank
not greater than one is its trace. The latter is ${\rm tr}\,{\bf P} = x^{1-k}-x$,
since ${\rm det}( x{\bf 1} + {\bf P}) = 1 $.)

On such leaves, the family of functions ${\rm tr}\,{\bf A}^n, n=1,\ldots , k-1$,
form a full set of commuting variables.
Let us introduce local coordinates on these symplectic leaves in the following
way.  Let $\lambda_1,\ldots,\lambda_k$ be the eigenvalues of the matrix ${\bf A}$
and $q_1,\ldots,q_k$ be the diagonal matrix elements of ${\bf B}$
in the basis in which ${\bf A}$ is diagonal. Imposing the condition
${\rm rk}\, {\bf P} \leq 1$
and conjugating ${\bf B}$ by a diagonal matrix one can make
${\bf B}$ take the form
$$
{\bf B}^i_j = \frac{\sqrt{q_i q_j}(1-x)}{ \lambda_i /
\lambda_j - x}~.\noeq{A}
$$
The functions $\lambda_i$ and $q_j$ are locally
 well defined functions on the symplectic leaves and the Poisson brackets
between them are

$$
\{\lambda_i \,,\lambda_j\} = 0,\noeq{A}\label{a1}
$$
$$
\{q_i\,,q_j\} = q_i q_j \,\frac{(\lambda_i + \lambda_j)}{(\lambda_i /
\lambda_j - x)(\lambda_j / \lambda_i - x)(\lambda_i - \lambda_j)} ~,
\noeq{A}\label{a2}
$$
$$
\{\lambda_i\,,q_j\} = \lambda_i q_j \delta_{i,j}~.\noeq{A}\label{a3}
$$

{\it Proof of the formulas eqs.(A\ref{a1}--\ref{a3}).} To simplify
calculations we assume for a moment that we are working with the group
$GL(k)$ rather than $SL(k)$. Having computed the Poisson structure on
the quotient space  by the lattice gauge group we can then restrict
it to the subspace corresponding
to $SL(k)$-connections since the latter space is a Poisson submanifold
in the whole quotient space.
The bivector defining the Poisson structure on ${\cal A}^l$ for the group
$GL(k)$ can be rewritten in the form
$$ B =
\frac{1}{2}\sum_{i,j,u,v \in \{1\ldots 4\}} E^{i (u)}_j \otimes E^{j (v)}_i
(\epsilon(u,v) + \epsilon(i,j)),\noeq{A} \label{Poiss}
$$ where
$\epsilon(i,j)$ is $-1, 0$, or $1$ if
$i$ is less, equal, or greater than $j$, respectively, and  $E^{i (u)}_j$ are
the standard $gl(k)$
generators acting on the $u$-th end of an edge.
(In our case  $E^{i (1)}_j$ acts on ${\bf A}$ from the left,  $E^{i (2)}_j$
acts on ${\bf B}$ from the left, $E^{i (3)}_j$ acts on ${\bf A}$ from the
right and $E^{i (4)}_j$ acts on ${\bf B}$ from the right.)

It is not practical, however, to compute the Poisson brackets between
$ \lambda_i$ and $ q_j$
using this bivector as it is, because it does not agree with the
diagonal form of the matrix ${\bf A}$. In order to avoid this
difficulty,
since we are interested only in computing Poisson brackets of gauge
invariant functions, we may change the bivector (A\ref{Poiss})
in such a way that it still defines the same Poisson bracket on the
space of gauge invariant functions. In other words, there are different
ways to write down Poisson brackets on the coset space
$G \times G / {\sl Ad}\,G$ in terms of a bivector on the space $G \times G$.
On the other hand, since we know that the projection $\pi : G \times G
\rightarrow G \times G / {\sl Ad}\,G$ is a Poisson map, i.e.
the bracket of gauge invariant functions is
gauge invariant, it suffices to compute the value of the bracket of
two such functions on a submanifold
$F \subset G \times G$ which
intersects each gauge orbit (i.e. each ${\sl Ad}\,G$-orbit)at least once.
In doing this way, we can simplify computations by
changing the bivector (A\ref{Poiss}) by terms vanishing on $F$.
As a prescription, one can formulate the following rule of allowed
modifications of the bivector.
One can add to any vector $E^{i(\alpha)}_j$ a vector
which is tangent to gauge orbits or vanishes on $F$.
The vectors tangent to the gauge orbits are just generators of the gauge
transformations, in our case $\sum_{u=1}^4 E^{i(u)}_j$. The vectors
vanishing on $F$ (which is in our case the space of connections with
${\bf A}$ diagonal) are, for example,
$\lambda_j E^{i(1)}_j + \lambda_i E^{i(3)}_j$.
Using these rules one can replace:
$$
E^{i (1)}_j \; \leadsto \; \frac{\lambda_i}{\lambda_j- \lambda_i}(E^{i
(2)}_j +E^{i (4)}_j ), \noeq{A}
$$
$$
E^{i (3)}_j  \; \leadsto \; \frac{\lambda_j}{\lambda_i-
\lambda_j}(E^{i (2)}_j +E^{i (4)}_j ). \noeq{A}
$$
By this trick the bivector $B$ can be transformed to the form

$$
B' = \sum_{i>j} E^{i (2)}_j \wedge E^{j(4)}_i
\frac{\lambda_i + \lambda_j}{2(\lambda_i-\lambda_j)}  +
\frac{1}{2}\sum_{i} E^{i (2)}_i \wedge E^{i(1)}_i + E^{i (3)}_i \wedge
E^{i(4)}_i ~. \noeq{A}
$$

Applying this bivector (which now leaves ${\bf A}$ diagonal) we get the
desired Poisson brackets. $\Box$

   The form of the brackets eqs.(\ref{a1}--\ref{a3}) is such that in order to
define the variables canonically conjugated to $\lambda_i$ we can just
multiply $q_i$ by  factors not depending on $q_i$. For example, one can take
the variables
$$
s_i = q_i \, x^{\frac{n-1}{2}}\left(\prod_{k,k \neq i}
\frac{(\lambda_k-\lambda_i)(\lambda_i-\lambda_k)}
{(\lambda_k-x\lambda_i)(\lambda_i-x\lambda_k)}\right)^\frac{1}{2}. \noeq{A}
$$

One can check by an explicit computation that these new variables, $s_i$, have
Poisson brackets
$$
 \{s_i\,,s_j\} =0, \noeq{A}
$$
$$
\{\lambda_i\,,s_j\} =
\lambda_i s_j \delta_{i,j}.\noeq{A}
$$

Using the formula
$$
\det {\bf B} = x^{\frac{n(n-1)}{2}}\prod_i q_i \prod_{i \neq
j}\frac{(\lambda_i-\lambda_j)}
 {(x\lambda_i-\lambda_j)} ~,\noeq{A}
$$
which can be easily proved by induction,  one can express the function $H =
{\rm tr}\, ({\bf B} + {\bf B}^{-1})$ in terms of  $\lambda_i$ and $s_i$;

$$
H = \sum_i(s_i+s_i^{-1} )\, x^{\frac{n-1}{2}}\left(\prod_{k,k \neq
i} \frac{(\lambda_k-\lambda_i)(\lambda_i-\lambda_k)}
{(\lambda_k-x\lambda_i)(\lambda_i-x\lambda_k)}\right)^\frac{1}{2}\noeq{A}
\label{ruij}
$$
that turns out to be just the Ruijsenaars hamiltonian.

 Note, that the Poisson structure on ${\cal A}^l = G \times G$ is nice
from various points of view. In particular, it is nondegenerate close
to the identity. The action of the group on this space by conjugation
is a Poisson one and has a well defined momentum map (in the sense of
Poisson-Lie groups) $\mu:G \times G \rightarrow G^*$
such that $\mu:({\bf A},{\bf B}) \mapsto
{\bf A}{\bf B}{\bf A}^{-1}{\bf B}^{-1}$.

\vspace{3mm}
{\large \bf Acknowledgements}\\
\noindent
This work was partially supported by the Grant No. RFBR-95-01-01101 and by
the Grant No. 96-15-96455 for the support of scientific schools.

\newpage
\epsfxsize16cm
\centerline{\epsfbox{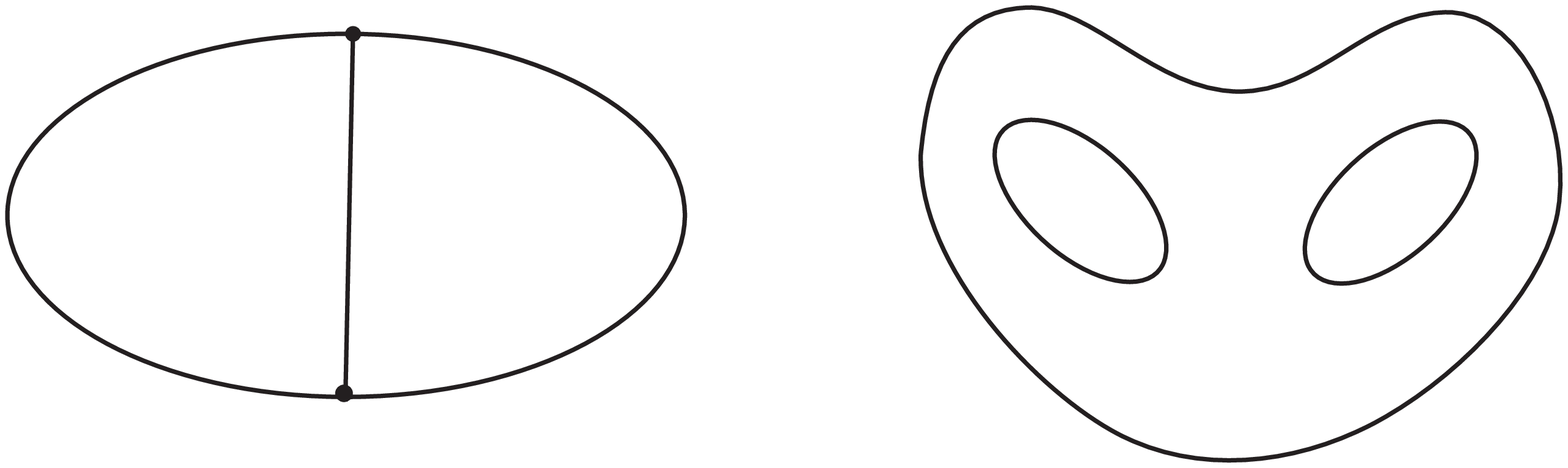}}
\vspace{0.5cm}
\centerline{ {\large\em a} \hspace{9cm} {\large\em b}}
\vspace{1cm}
\epsfxsize16cm
\centerline{\epsfbox{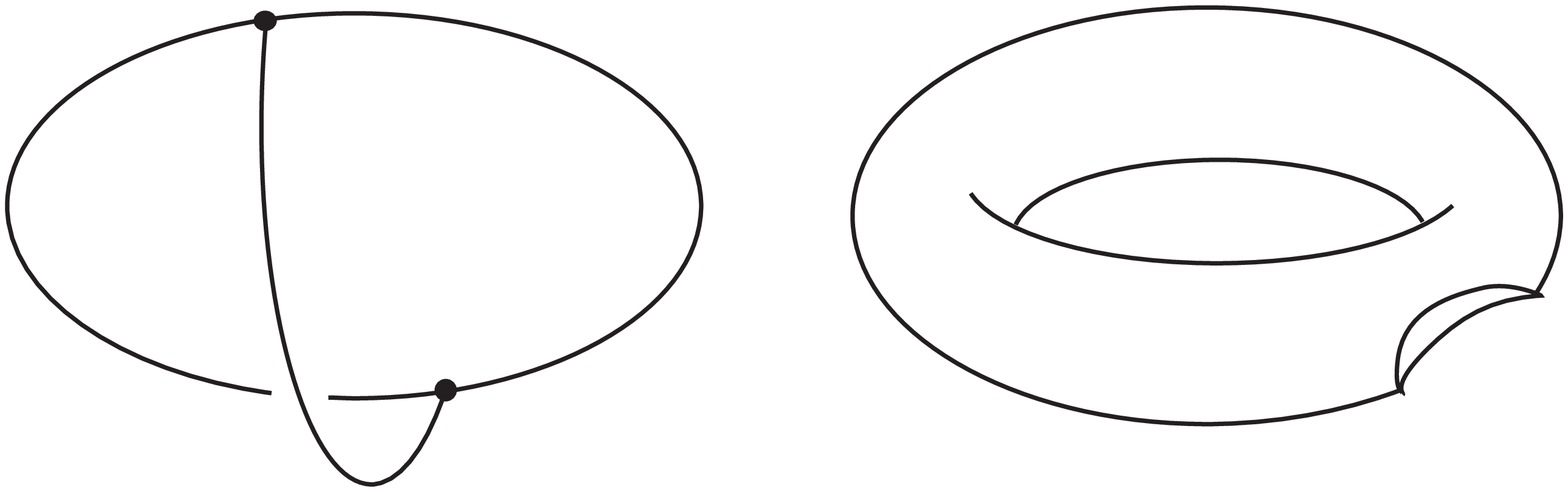}}
\vspace{0.5cm}
\centerline{ {\large\em c} \hspace{9cm} {\large\em d}}
\vspace{1cm}
\centerline{Fig. 1}
\vspace{0.25cm}
\begin{minipage}{15cm}
\centerline
{Examples of fat graphs and surfaces corresponding to them.}

\centerline
{The cyclic orders at vertices are understood to be counterclockwise.}

\centerline
{The graph {\large\em (a)} gives a disk with two holes {\large\em (b)}.}

\centerline
{The graph {\large\em (c)} gives a torus with one hole {\large\em (d)}.}
\end{minipage}

\newpage
\epsfxsize16cm
\centerline{\epsfbox{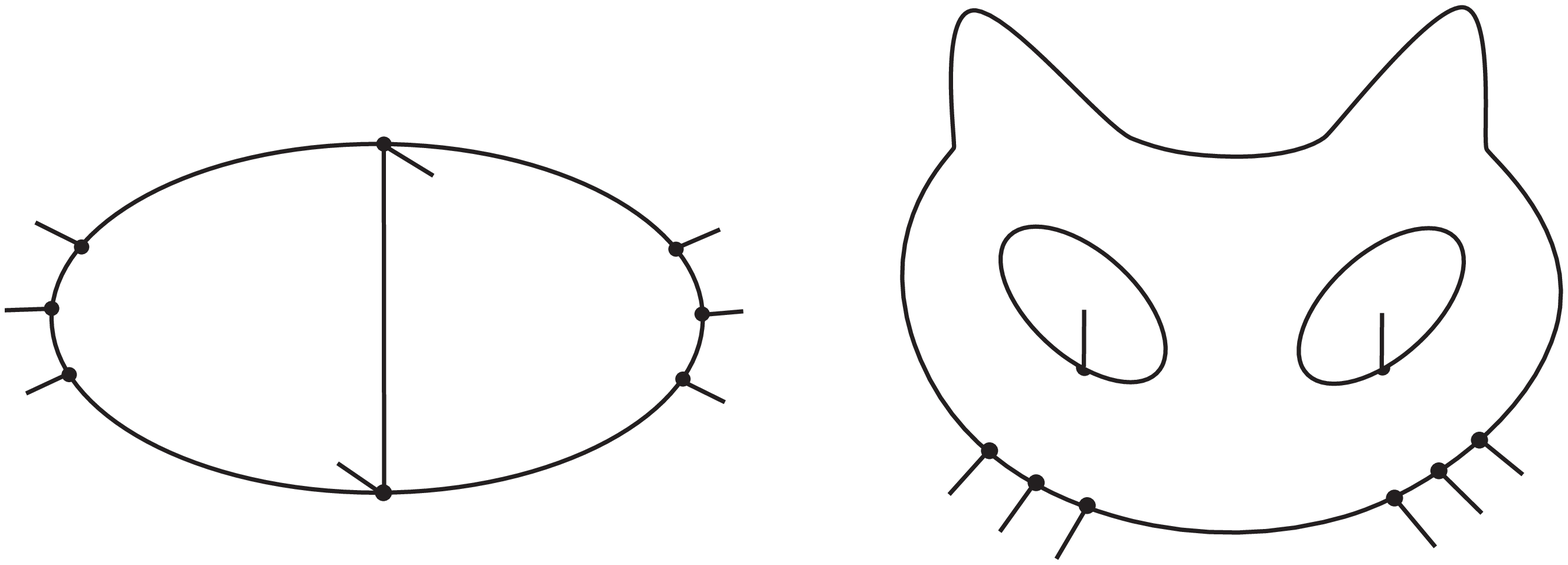}}
\vspace{0.5cm}
\centerline{ {\large\em a} \hspace{9cm} {\large\em b}}
\vspace{1cm}
\epsfxsize16cm
\centerline{\epsfbox{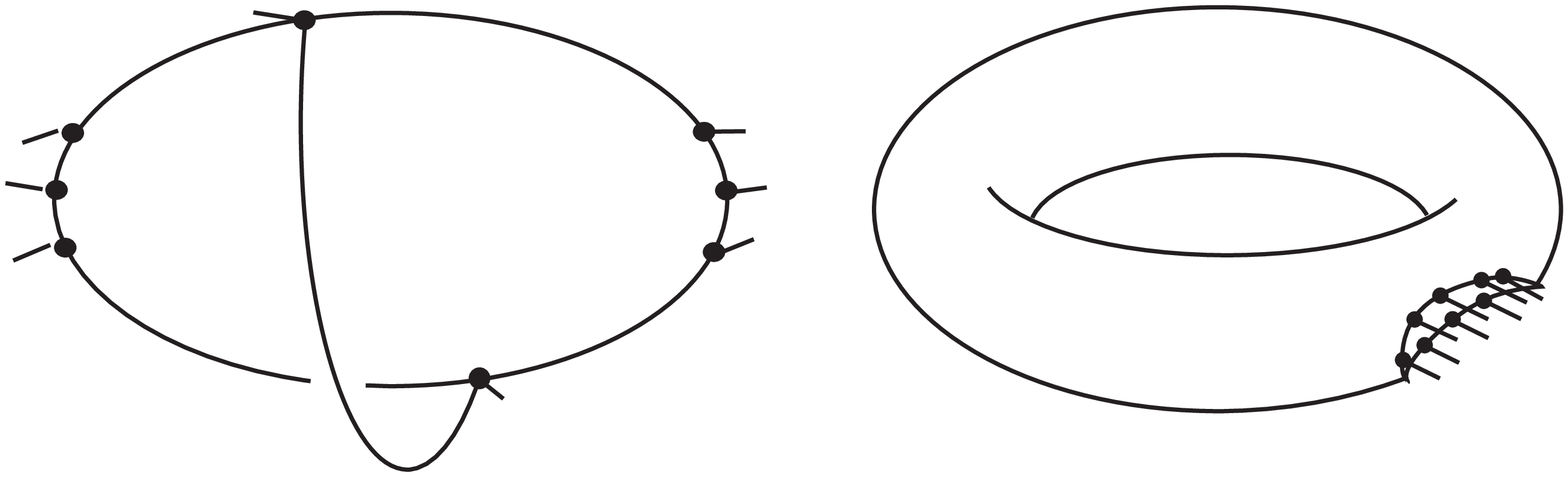}}
\vspace{0.5cm}
\centerline{ {\large\em c} \hspace{9cm} {\large\em d}}
\vspace{1cm}
\centerline{Fig. 2}
\vspace{0.25cm}
\begin{minipage}{15cm}
\centerline
{Examples of ciliated fat graphs and corresponding ciliated surfaces.}

\centerline
{Cilia are indicated by small strokes at the vertices.}

\centerline
{The graph {\large\em (a)} gives a disk with two holes {\large\em (b)}.}

\centerline
{The graph {\large\em (c)} gives a torus with one hole {\large\em (d)}.}
\end{minipage}

\newpage
\epsfxsize16cm
\centerline{\epsfbox{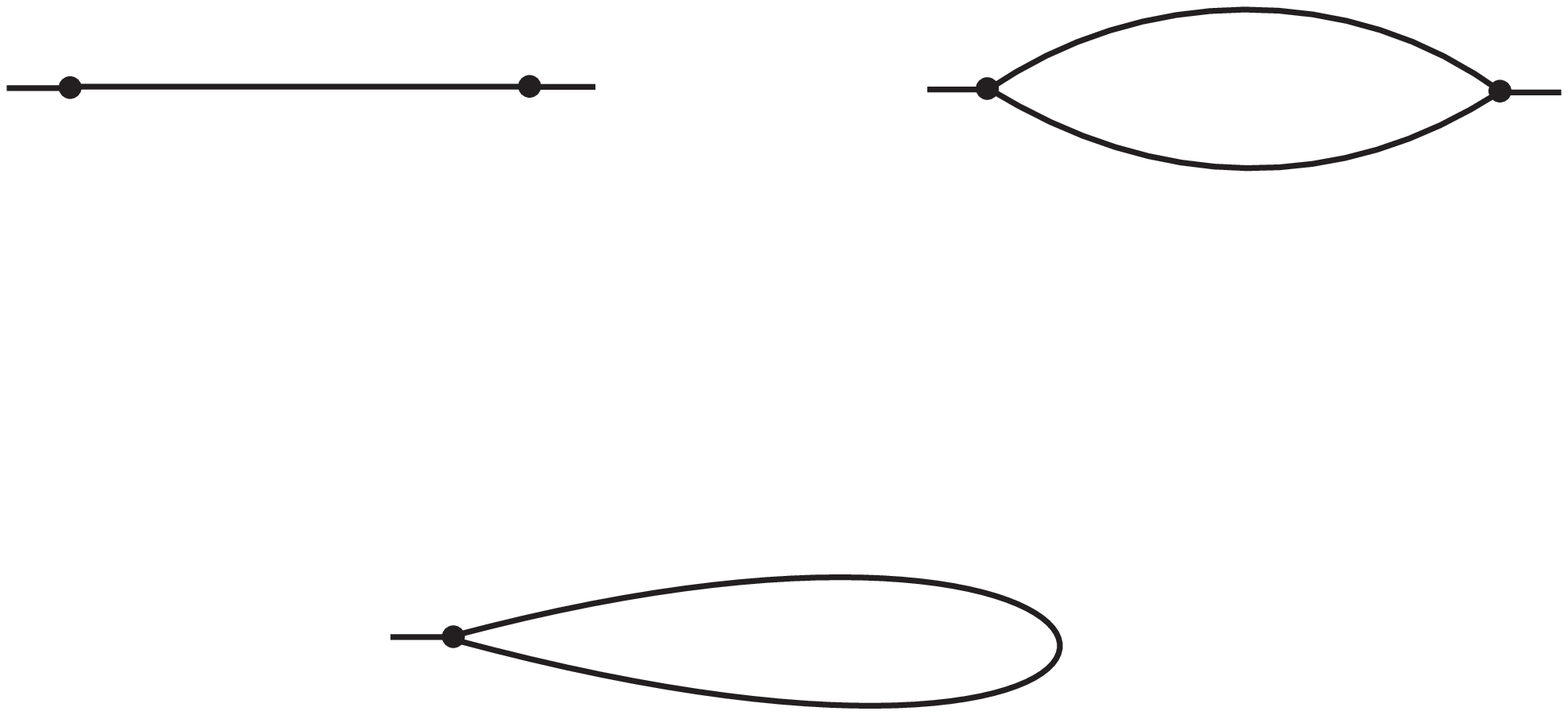}}
\vspace{0.5cm}
\begin{picture}(0,0)(0,0)
\unitlength1cm
\put(2.3,6.3){\makebox(0,0)[cc]{\large\em a}}
\put(7.5,0.5){\makebox(0,0)[cc]{\large\em c}}
\put(12.3,6.3){\makebox(0,0)[cc]{\large\em b}}
\end{picture}

\vspace{0.75cm}
\centerline{Fig. 3}
\vspace{0.25cm}

\begin{minipage}{15cm}
\centerline{The graphs corresponding to}

~~~~~~~~~~~~~~~~~~~~~~~~~~~~~~~~~~
{\large\em (a)} the Poisson-Lie group $G$,

~~~~~~~~~~~~~~~~~~~~~~~~~~~~~~~~~~
{\large\em (b)} its double $D \simeq G \times G$,

~~~~~~~~~~~~~~~~~~~~~~~~~~~~~~~~~~
{\large\em (c)} its dual Poisson-Lie group $G^*$.
\end{minipage}

\vspace{2cm}

\epsfxsize16cm
\centerline{\epsfbox{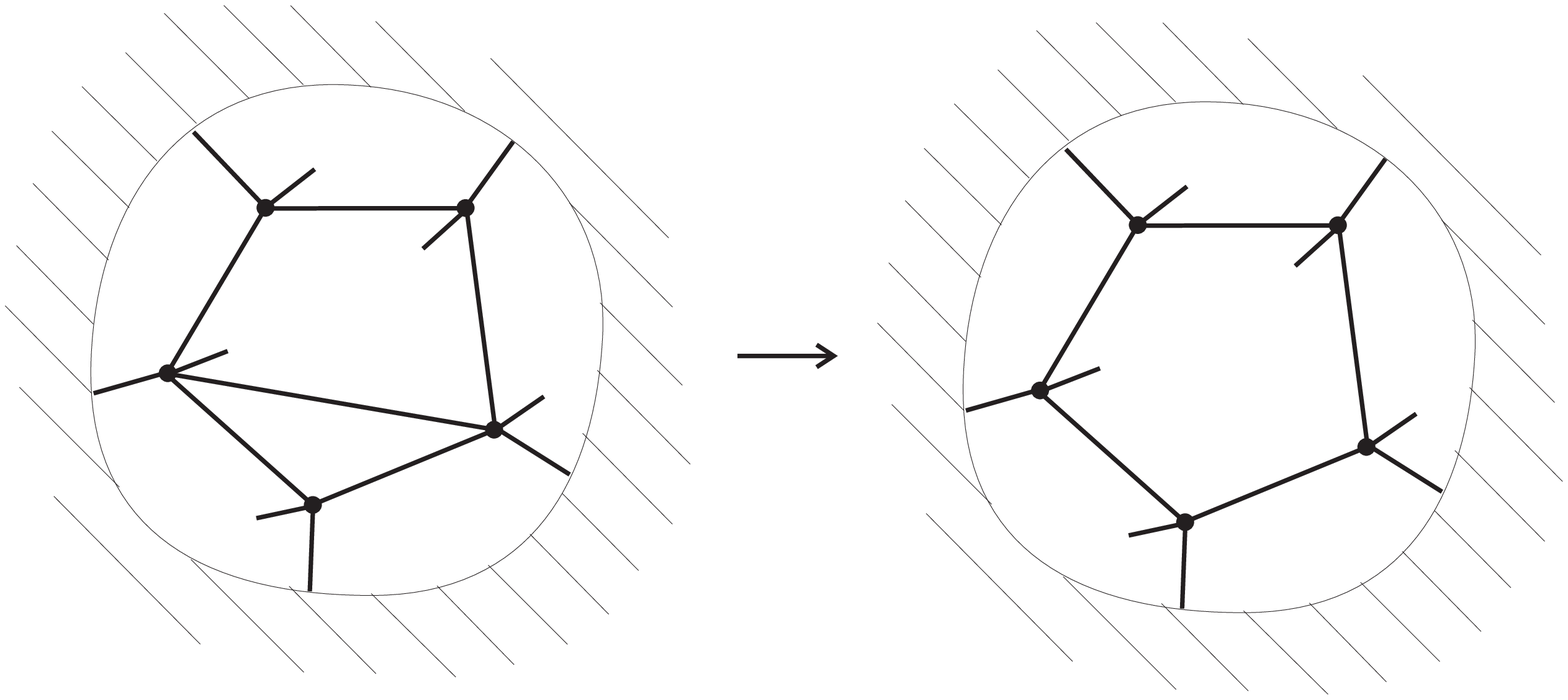}}
\centerline{Fig. 4}
\vspace{0.25cm}
\begin{minipage}{15cm}
\centerline{Operation of erasing an edge.}
\centerline
{The shaded region represents the remainder of the graph.}
\end{minipage}

\newpage
\epsfxsize16cm
\centerline{\epsfbox{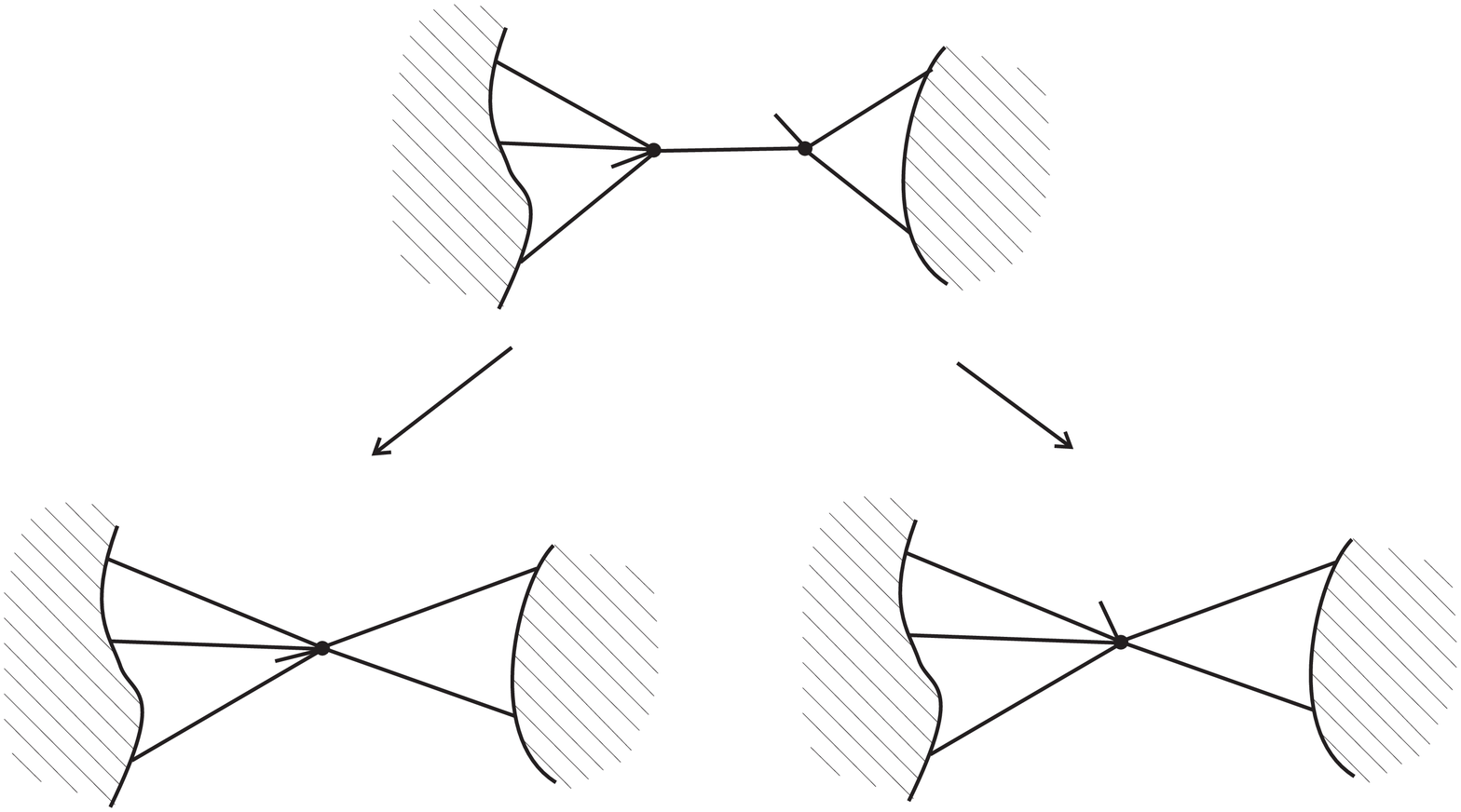}}
\begin{picture}(0,0)(0,0)
\unitlength1cm
\put(4.6,4.8){\makebox(0,0)[cc]{$L$}}
\put(9.5,4.8){\makebox(0,0)[cc]{$R$}}
\put(6.7,7.4){\makebox(0,0)[cc]{$n_L$}}
\put(8,7.4){\makebox(0,0)[cc]{$n_R$}}
\end{picture}
\vspace{0.5cm}

\centerline{Fig. 5}

\begin{center}
Operations of contractions of an edge.

$L$ and $R$ are the two different ways of contraction.

$L$ corresponds to factoring by gauge transformations at the  vertex $n_R$.

$R$ corresponds to factoring by gauge transformations at the  vertex $n_L$.
\end{center}

\vspace{2.5cm}

\epsfxsize16cm
\centerline{\epsfbox{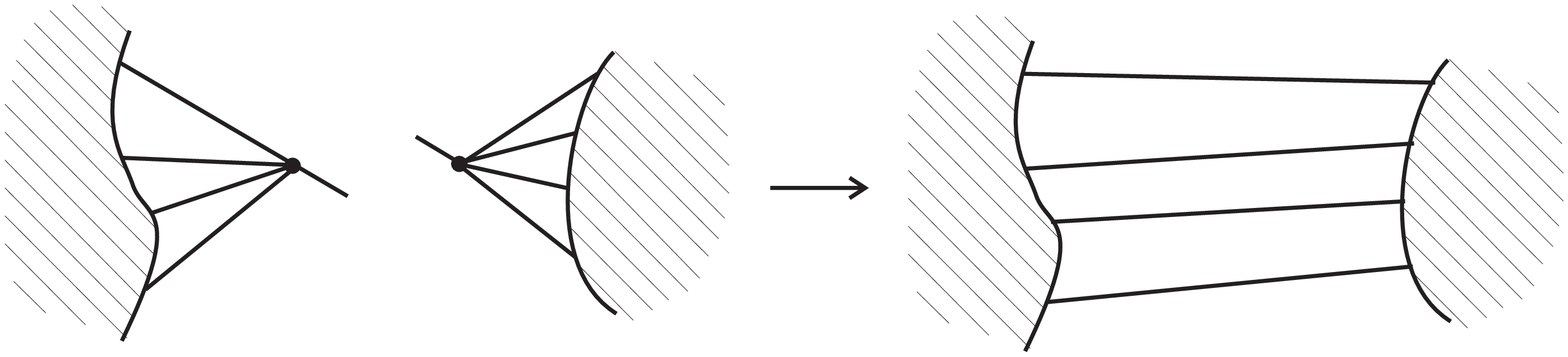}}
\vspace{1cm}
\centerline{Fig. 6}
\vspace{0.25cm}
\centerline{Operation of gluing graphs.}

\newpage
\epsfxsize16cm
\centerline{\epsfbox{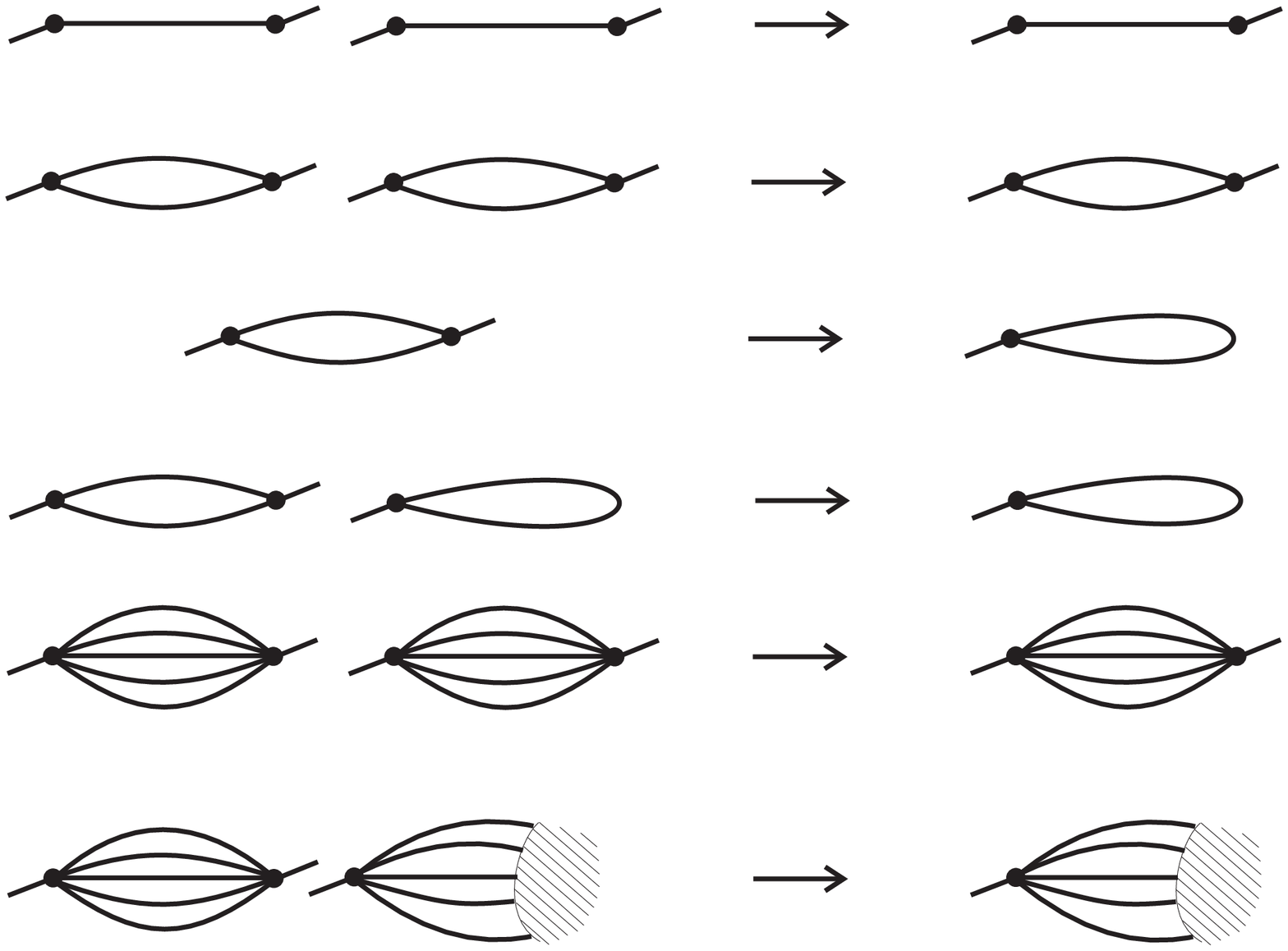}}
\unitlength1cm
\begin{picture}(0,0)(0,0)
\put(-1,12){\makebox(0,0)[cc]{\large\em a}}
\put(-1,10){\makebox(0,0)[cc]{\large\em b}}
\put(-1,8.1){\makebox(0,0)[cc]{\large\em c}}
\put(-1,6.1){\makebox(0,0)[cc]{\large\em d}}
\put(-1,4.1){\makebox(0,0)[cc]{\large\em e}}
\put(-1,1.3){\makebox(0,0)[cc]{\large\em f}}
\end{picture}

\vspace{2cm}
\centerline{Fig. 7}

\begin{center}
Some particular cases of gluing graphs which correspond
to natural operations in Poisson-Lie groups:
\end{center}

\begin{minipage}{15cm}
~~~~~~~~~~~~~~~~~~~~~~~~~~~~~~~~~~~~~
{\large\em (a)} multiplication in $G$,

~~~~~~~~~~~~~~~~~~~~~~~~~~~~~~~~~~~~~
{\large\em (b)} multiplication in $D$,

~~~~~~~~~~~~~~~~~~~~~~~~~~~~~~~~~~~~~
{\large\em (c)} projection $D \rightarrow G^*$,

~~~~~~~~~~~~~~~~~~~~~~~~~~~~~~~~~~~~~
{\large\em (d)} action of $D$ on $G^*$,

~~~~~~~~~~~~~~~~~~~~~~~~~~~~~~~~~~~~~
{\large\em (e)} multiplication in the 5-uble,

~~~~~~~~~~~~~~~~~~~~~~~~~~~~~~~~~~~~~
{\large\em (f)} action of the 5-uble on a space

~~~~~~~~~~~~~~~~~~~~~~~~~~~~~~~~~~~~~~~~~~~
of graph connections.
\end{minipage}

\newpage
\epsfxsize16cm
\centerline{\epsfbox{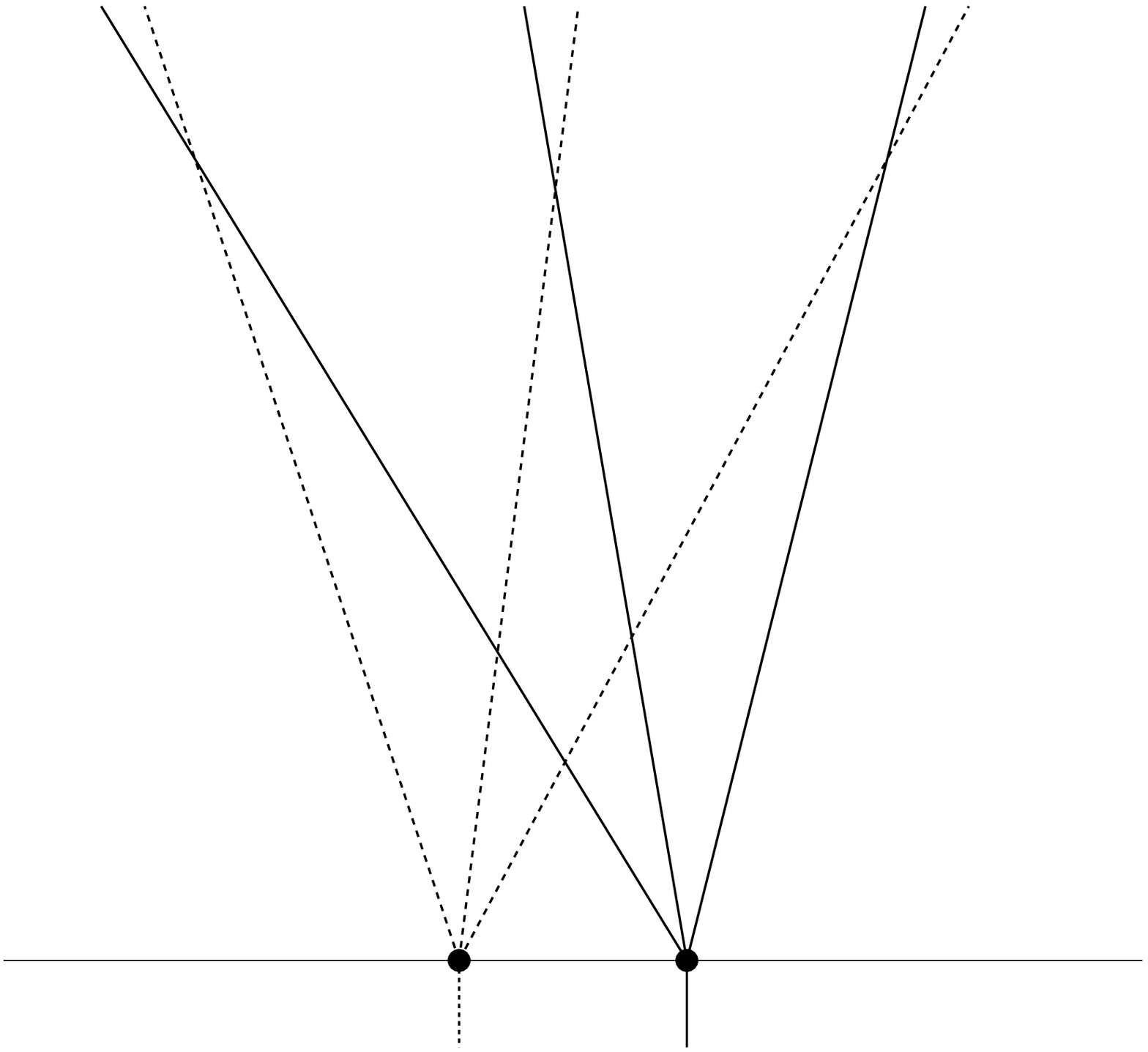}}
\vspace{2cm}
\centerline{Fig. 8}
\begin{center}
The special way of deforming a graph drawn on a surface which
gives two transversal graphs; the deformed graph is shown by  the
broken line.
\end{center}

\newpage
\epsfxsize14cm
\centerline{\epsfbox{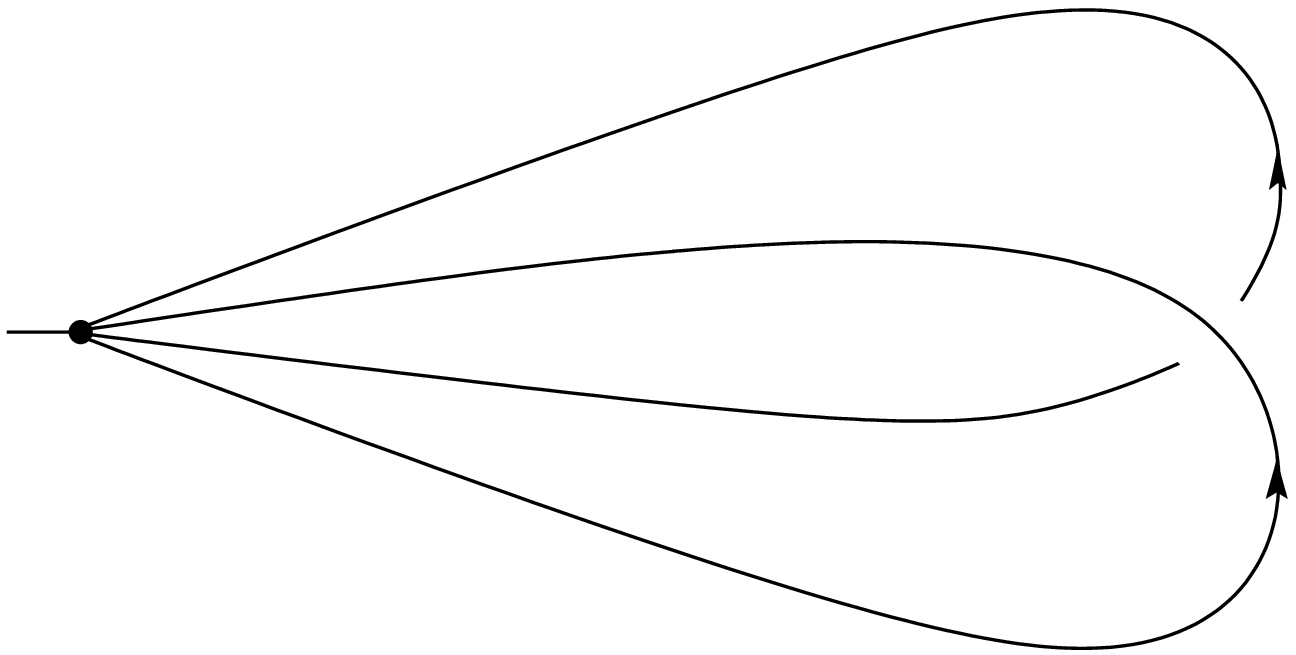}}
\vspace{1cm}

\begin{picture}(0,0)(0,0)
\unitlength1cm
\put(5,6.8){\makebox(0,0)[cc]{{\sl 1}}}
\put(5,5.9){\makebox(0,0)[cc]{{\sl 2}}}
\put(5,4.2){\makebox(0,0)[cc]{{\sl 3}}}
\put(5,3.2){\makebox(0,0)[cc]{{\sl 4}}}
\put(14.7,6.8){\makebox(0,0)[cc]{$A$}}
\put(14.7,3.5){\makebox(0,0)[cc]{$B$}}
\end{picture}

\centerline{Fig. 9}
\vspace{0.25cm}
\centerline{The ciliated graph corresponding to a holed torus.}

\vspace{3cm}

\epsfxsize14cm
\centerline{\epsfbox{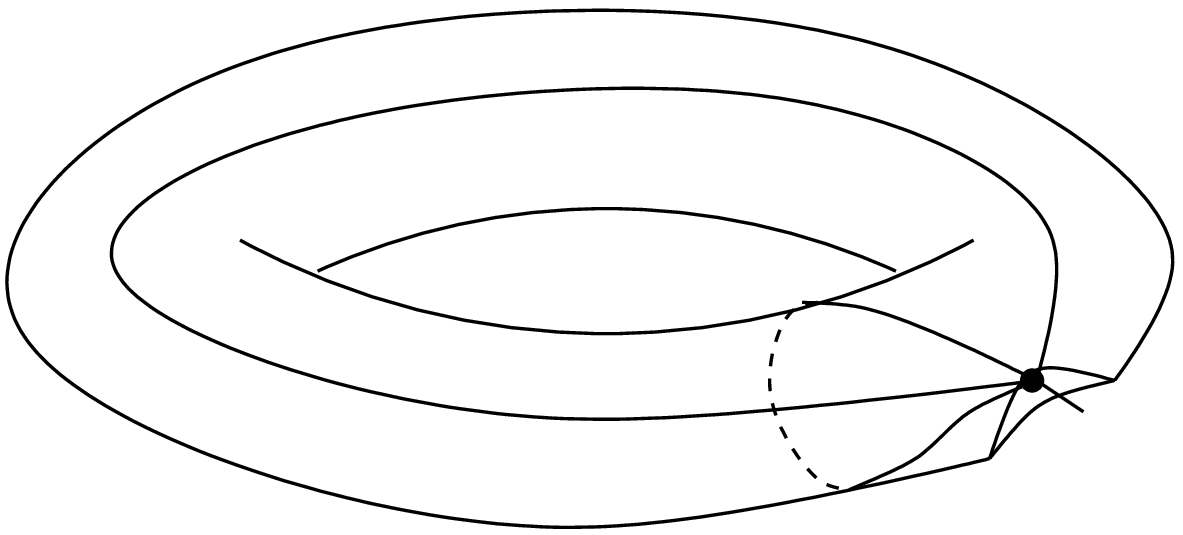}}
\vspace{0.25cm}
\centerline{Fig. 10}
\vspace{0.25cm}
\centerline{The same graph drawn on the holed torus.}
\begin{picture}(0,0)(0,0)
\unitlength1cm
\put(3.3,3.8){\makebox(0,0)[cc]{$A$}}
\put(11.3,4.8){\makebox(0,0)[cc]{$B$}}
\end{picture}

\end{document}